\documentclass[onefignum,onetabnum]{siamart250211}

\usepackage{amsfonts}
\usepackage{subcaption}
\usepackage{graphicx}
\usepackage{epstopdf}
\usepackage{algorithmic}
\ifpdf
  \DeclareGraphicsExtensions{.eps,.pdf,.png,.jpg}
\else
  \DeclareGraphicsExtensions{.eps}
\fi
\usepackage{enumerate}

\newcommand{\R}{\mathbb{R}}

\newcommand{\N}{\mathbb{N}}

\newcommand{\laplace}{\Delta}

\DeclareMathOperator{\trace}{trace}
\DeclareMathOperator{\hessian}{H}
\DeclareMathOperator{\diver}{div}

\DeclareMathOperator{\supp}{supp}

\newcommand{\Lan}{{\mathcal L}} 

\newcommand{\Ord}{{\mathcal O}}

\newsiamremark{remark}{Remark}
\newsiamremark{hypothesis}{Hypothesis}
\crefname{hypothesis}{Hypothesis}{Hypotheses}
\newsiamthm{claim}{Claim}

\headers{Harmonic Eigenfunctions for MPI Reconstruction}
{T. März, V. Gapyak, and A. Weinmann}

\title{
Higher Order Regularization using Harmonic Eigenfunctions for Model-Based Reconstruction in Magnetic Particle Imaging
}

\author{
Thomas März
\thanks{
ACIDA~Lab,~Hochschule~Darmstadt;~(\email{thomas.maerz@h-da.de}, \email{vladyslav.gapyak@h-da.de}, \email{andreas.weinmann@h-da.de}).
}
\and 
Vladyslav Gapyak
\and 
Andreas Weinmann 
}

\usepackage{amsopn}

\begin{document}

\maketitle

\begin{abstract}
Magnetic Particle Imaging (MPI) is a recent imaging modality 
where superparamagnetic nanoparticles are employed as tracers. 
The reconstruction task is to obtain the spatial particle distribution from
a voltage signal induced by the particles.
Generally, in computational imaging variational reconstruction techniques are 
common and rely on a mathematical model to describe the underlying physics.
For the MPI reconstruction task we propose a model-based variational 
reconstruction technique which incorporates a higher order regularizer, 
where the regularizer is diagonalized by harmonic eigenfunctions. 
The proposed image reconstruction algorithm features two major stages:
in the first stage, the core stage, the components of the MPI core response 
are reconstructed. This is the MPI-specific data approximation task 
which we formulate as a variational problem
incorporating the higher order regularizer.
The relationship between the particle distribution, the MPI core response 
and the measured data is given by a mathematical model which was introduced 
in our earlier research. According to this model the MPI core response is 
tied to the particle distribution by convolution.
Therefore the outcome of the core stage yields the data for the second stage, 
the deconvolution stage, in which the final reconstructed image is produced 
by solving an ill-posed deconvolution problem 
in a robust way relying on earlier research.
Interestingly, the quality of the final image depends significantly on the 
quality of the result of the core stage.
A contribution is thus the enhancement of the core stage 
via higher order regularization. 
We provide a theoretical foundation for our approach and 
demonstrate its benefit with numerical examples.
\end{abstract}

\begin{keywords}
Image Reconstruction, Data Approximation, Energy Minimization, Regularization,
Laplacian, Bi-Laplacian, Harmonic Eigenfunctions
\end{keywords}

\begin{MSCcodes}
65K10,65R32,65T40,92C55
\end{MSCcodes}

\section{Introduction}

Magnetic Particle Imaging (MPI) is a tracer-based imaging modality that utilizes 
superparamagnetic nanoparticles as tracer material.
The non-linear magnetization response of such particles, when dynamic magnetic 
fields are applied, is harnessed to encode
an unknown spatial particle distribution into a time-dependent voltage signal. 
During the scanning procedure the nanoparticles, which were injected 
in the specimen being scanned, 
induce a voltage in the receive coils of the MPI-scanner. 
The received voltage signal is recorded as a time series and constitutes 
the data from which the 
spatial distribution of the nanoparticles is to be reconstructed. 
The spatial particle distribution is the image which provides a means to 
inspect or visualize interior structures of the specimen, e.g., 
the structure of blood vessels or the location of tumors. 
Since the introduction of MPI by Gleich and Weizenecker~\cite{gleich2005original}
in 2005 many different medical applications have been proposed and investigated, 
among which are cancer detection~\cite{Song2018} 
and cancer imaging~\cite{Tay2021},
blood flow imaging~\cite{Franke2020BloodFlow}, or 
stem cell tracing~\cite{connell2015advancedcellTherapies}. 
Further applications of MPI can be found in~\cite{Yang2022}. 

Regarding the reconstruction task there are essentially two type of approaches:
there are model-based approaches, relying on mathematical models of 
the underlying physics to emulate the MPI imaging operator, 
and measurement-based approaches, employing a system matrix acquired 
via calibration.
The calibration procedure in the measurement-based approach is performed 
by scanning delta-concentrations 
in known positions on a voxel grid. The measured scan data, i.e., 
the measured voltage signals of these delta-concentrations, are collected 
into a system matrix.
Given the scan data of a general specimen, the respective particle 
distribution is reconstructed 
from a system of linear equations by regularized inversion of the 
system matrix~\cite{Rahmer_etal2012}.
Because said linear system is ill-conditioned, regularization is necessary.
Among the techniques for regularized inversion there are classical techniques 
such as Tikhonov regularization~\cite{Weizenecker_etal2009,Knopp_etal2010ec} 
and $\ell^1$-priors~\cite{storath2016edge},
often combined with the Kaczmarz method, as well as 
machine learning-based approaches
such as Deep Equilibrium Models~\cite{gungor2023deqmpi} or 
the Deep Image Prior (DIP)~\cite{dittmer2020deep,Yin2023dipSM}.   
Recently, so-called Plug-and-Play (PnP) priors as suggested 
by~\cite{venkatakrishnan2013pnp}
have been successfully adopted for the MPI-reconstruction task 
\cite{askin2022pnp,gapyak2025ell,gapyak2025fast}
and employ machine learning-based denoisers in lieu of classical 
regularization steps 
within iterative reconstruction schemes.
Because the calibration procedure of the measurement-based approach is 
time-consuming, storage-intensive, and a priori depending
on the grid resolution (current grids in 3D are about $20 \times 20 \times 20$ 
voxels \cite{knopp2020openmpidata}),
model-based approaches have been considered early on, e.g., 
in~\cite{rahmer2009signal,knopp2009model,knopp20102d}. 
While \cite{rahmer2009signal,knopp2009model} focus on the 1D scenario in a 
Field Free Point (FFP) setup, the later X-Space formulation
of~\cite{GoodwillConolly2010,GoodwillConolly2011,Goodwill2012Proj} 
tackles also the 2D and 3D scenarios.
In~\cite{knopp20102d} the 2D MPI problem is approached by simulating 
the system matrix.
These approaches are all based on the Langevin theory of 
paramagnetism~\cite{jackson_classical_1999}.
Later contributions include more details about particle 
features~\cite{kluth2018mathematical,maass2024equilibrium} 
or provide more detailed models of the applied magnetic 
fields~\cite{bringout2020new}.
More recently, it has been demonstrated in \cite{maass2024magnetic} 
that for immobilized particles the Lagevin model is sufficient to model the magnetization 
behavior of the paramagnetic nanoparticles.
 
Inspired by the X-space formulation the present authors in~\cite{marz2016model} 
decomposed the MPI signal encoding into two consecutive steps.
From this decomposition reconstruction formulae in 2D and 3D were 
derived and a first two-stage reconstruction algorithm was distilled. 
The two-stage strategy was further developed by the present authors to provide 
various points of flexibility for enhancement 
of the reconstruction quality~\cite{gapyak2022mdpi,marz2022amee}, 
to include multi-patch scan data~\cite{gapyak2023multipatch}, 
the reconstruction from partial data~\cite{marz2022icnaam,gapyak2023amee},
and the Field Free Line (FFL) MPI-setup \cite{gapyak2025ffl3d}.
The flexibility to handle different scans of the same object is due to a 
variational formulation~\cite{gapyak2022mdpi,marz2022amee}
of the core stage, in which the model of the signal encoding is employed 
in the data fidelity term to reconstruct the MPI core response.
In the deconvolution stage the actual particle distribution is reconstructed 
from the result of the core stage by regularized deconvolution 
employing the convolution kernel derived in~\cite{marz2016model}.
Enhanced reconstruction quality was achieved~\cite{gapyak2023multipatch} 
employing a combination of TV-regularization and a sparsity-supporting 
prior as well as positivity constraints.
More recently, the present authors successfully applied the model-based 
two-stage approach to reconstruct images from real 
scanning data \cite{gapyak2025fast}.

In this article we focus on the MPI-specific core stage of the algorithm to 
contribute further enhancement of the overall reconstruction quality by employing 
higher order regularization in the core stage of the algorithm.
The quality of the final image depends strongly on the quality of the outcome 
of the core stage, which serves as data for the deconvolution stage.
The task of the core stage is to reconstruct the MPI core response, 
which is the response $A[\rho]$ of the MPI core operator $A$ to the 
yet unknown particle distribution $\rho$. 
The MPI core response is a smooth function of type
$A[\rho]:\Omega \to \R^{n \times n}$, where $\Omega = [-1,1]^n$ is the 
spatial domain and $n \in \{2,3\}$ for 2D or 3D MPI.
The collected data $s(t_i)$ is a time series and, in an ideal scenario,
is related to  the MPI core response by $s(t) = A[\rho](r(t)) v(t)$, 
where $r(t)$ is the motion of the FFP, i.e. the scanning curve, 
and $v(t)$ its velocity.
The reconstruction of $A[\rho]$ is a data approximation problem, where the 
data fidelity term models the relationship $s(t_i) = A[\rho](r(t_i)) v(t_i)$
for the discrete time series.
We consider FFP motions, which are employed in practice 
(e.g. Open MPI Data \cite{knopp2020openmpidata})
and are usually Lissajous curves.
A specific feature of this FFP motion is that big portions of the 
field of view (FoV), i.e. the domain $\Omega$, are never visited by 
the FFP as shown Fig.~\ref{fig:liss:scans} (left). 
Because of this the distribution of data sample locations
is very sparse, which is a challenge for reconstruction. 
In addition, increasing the data sampling rate does not mitigate this 
challenge since sample locations are confined to lie on the scanning 
trajectory, which stays unchanged.
Because of the data sparsity, we need regularization, which is able 
to fill the gaps appropriately and to handle the noise in the data. 

In the preliminary works \cite{gapyak2022mdpi,marz2022amee,gapyak2023multipatch} 
we employed first order regularization, namely the squared $L^2$-norm of 
the gradient. Here, we noticed that the reconstruction of $\trace(A[\rho])$ 
in 2D looks very spiky like a circus tent-shaped interpolation surface 
with ``tent poles" near the data sampling locations $r(t_i)$.
This effect is likely due to the maximum principle of the Laplacian 
\cite{evans2022partial} which appears in the Euler-Langrange equation.
This effect can be mitigated by collecting additional data \cite{gapyak2022mdpi},
but in a sparse data scenarios we need more appropriate regularization.

Inspired by the multivariate smoothing spline 
\cite{cox1984multivariate,hastie2009elements,wang2011smoothing}
we employ higher order regularization. In particular, 
we consider second order regularization, 
namely the squared $L^2$-norm of the Laplacian.
This regularizer appears also as the energy term in the theory of 
thin plate splines \cite{friedrichs1928randwert,reddy2006theory,timoshenko1959theory}.
In the setup here, we work on an axis-parallel rectangular 
finite domain $\Omega$ and employ zero Neumann boundary conditions 
which are typically employed in image and geometry 
processing \cite{stein2018natural}.
Because of the particular domain, it is quite natural to expand 
solutions in terms of eigenfunctions
of the underlying operators, which are the Laplacian (first order) 
and the Bi-Laplacian (second order)
supplemented with boundary conditions.
Depending on the choice of boundary conditions, there exist eigenfunctions 
with separated variables,
either tensor-product sine functions or tensor-product cosine functions.
In \cite{courant2004methods} the derivation of the eigenfunctions 
for the Laplacian with 
either zero Dirchlet or zero Neumann boundary conditions is given.
For a simply supported thin plate the tensor-product sine functions are used
to construct the Navier solution \cite{reddy2006theory}.
In the case considered here the tensor-product cosine functions turn 
out as harmonic eigenfunctions
to expand the solution.
 
\section*{Contributions}
In this paper we propose a higher order regularization technique to enhance
the overall reconstruction quality of model-based reconstruction in MPI.
The basic idea was shortly communicated in the conference proceeding \cite{marz2023icnaam}, 
but without details on the method or its analysis, which we provide now.
More precisely, we make the following contributions: 
\begin{enumerate}[(i)]
\item 
We integrate Laplacian regularization, the $L^2$-norm of the Laplacian,
in the core stage of our two-stage algorithm for model-based MPI reconstruction. 
The second order Laplacian regularization yields a smoother reconstruction 
of the MPI core response, which is an 
analytic function according to the theory \cite{marz2016model}.
This helps significantly to improve the overall reconstruction results.
\item 
We contribute a novel two-stage algorithm, where the design of the core stage
uses tensor-product cosine functions opposed to finite difference schemes.
These cosine functions are the eigenfunctions of the Bi-Laplacian on a 
square with enforced zero Neumann boundary conditions which is 
the underlying operator of the employed Laplacian regularizer.
This results in a simpler description of the regularizer.
Moreover, it improves the performance of the method on modern GPUs. 
\item
We provide a theoretical foundation of our approach.
Laplacian as well as Hessian regularizers (which are common in image and geometry processing)
both feature the Bi-Laplacian as underlying operator,
but the choice of boundary conditions is more intricate and affects the eigenfunctions.
We discuss the existence of tensor-product eigenfunctions of the 
Bi-Laplacian on a square and reveal the importance of the chosen boundary conditions.
To the best of our knowledge, this has not been covered by the literature.  
In particular, we show that in the case of all natural boundary conditions only for the eigenvalue $\mu=0$ 
tensor-product eigenfunctions are available (Thm.~\ref{theo:EigBiLapAllNatSep}). 
Moreover, in this case the underlying operator has an infinite dimensional kernel (all harmonic functions,
Thm.~\ref{theo:BiLapAllNatKernel}).
When enforcing at least partially Neumann-zero or Dirchlet-zero boundary conditions,
we obtain either cosine or sine tensor-product eigenfunctions.
While with all enforced boundary conditions, there a no tensor-product eigenfunctions 
(Thm.~\ref{theo:BiLapEigenAllEnf}) at all.
Finally, we show, that in our setup Laplacian and Hessian regularization are actually identical 
(Thm.~\ref{theo:LapReguEquivHessRegu} and Thm.~\ref{theo:R2R3sameEig}).
\item 
We demonstrate the potential of our method by numerical experiments. 
\end{enumerate}

\section*{Outline}
In Sect.~\ref{sect:Model} a concise review 
of the mathematical model of MPI and the reconstruction formulae are given. 
In Sect.~\ref{sect:MethodAlgo} we describe the methodology and the algorithms.
In particular, we formulate the core stage of our two-stage algorithm in terms of harmonic eigenfunctions
and contrast this version of the preliminary algorithm with the novel 2nd order algorithm.
Sect.~\ref{sect:Theory} provides the theoretical foundations:
we discuss the existence of tensor-product eigenfunctions of the 
Bi-Laplacian given different choices of boundary conditions.
Moreover, we prove that the energies given by the squared $L^2$-norm of 
the Laplacian
and given by the squared $L^2$-norm of the Hessian are identical in 
the case of a rectangular domain
with enforced zero Neumann boundary conditions.  
In Sect.~\ref{sect:Experiments} we discuss the results of the 
computational experiments. 
We show with different experiments the gain in reconstruction quality.
Finally, Sect.~\ref{sect:Conclusion} presents our conclusions.

\section{Mathematical Model and Core Components of the Method}\label{sect:Model}

\subsection{Signal Encoding}
In this article
we consider the field free point (FFP) setup of an MPI-scanner, 
which is very common \cite{knopp2020openmpidata}.
In this setup an applied dynamic magnetic field steers the FFP $r(t)$ to 
scan a specimen which has been prepared
with superparamagnetic particles and placed in the scanner.
The FFP triggers locally the non-linear magnetization response of 
the superparamagnetic particles.
Using the Langevin theory \cite{chikazumi1978physics,jiles1998introduction}
to model the magnetization behavior of the particles and
combining it with Faraday's law, we obtain
a basic model \cite{BuzugKnopp2012,GoodwillConolly2010,GoodwillConolly2011}
to describe the part of the received signal $s(t)$ (the voltage induced 
in the recording coils)
which depends on the spatial distribution of the particles.
We work with a mathematically transformed version \cite{marz2016model} 
of the basic model which reads as
\begin{align}\label{eqn:signalTrafo}
s(t) &= \frac{d}{dt} \int\limits_{\R^n} \rho(x) \; \Lan \left( \frac{|r(t)-x|}{h} \right) \; \frac{r(t)-x}{|r(t)-x|} \; dx.
\end{align}
In the considered setups the dimension $n$ is either 2 or 3.
The data $s(t)$ depends linearly on the underlying signal $\rho$, 
which is the spatial
distribution of the particles and needs to be reconstructed from the data.
The particle distribution $\rho$ is confined to the so-called field of view
$\Omega = [-1;1]^n$,
i.e., the support $\supp \rho = \overline{\{ x : \rho(x) \neq 0\}}$ of 
$\rho$ is contained in $\Omega$.
In Eq.~\eqref{eqn:signalTrafo} the Lagevin-function 
$\Lan(\xi) = \coth(\xi) - 1/\xi$
describes the generic magnetization behavior.
It depends also on a dimensionless resolution parameter $h > 0$,
where smaller values of $h$ mean less blur.
The parameter $h$ is a combination of different physical parameters
(cf. \cite{marz2016model}) which correspond to properties of the
particles (diameter, saturation magnetization), the applied field 
(strength of the gradient) and the scanner (size of the field of view).

\subsection{Core Operator}\label{sect:CoreOp}
The MPI core operator acting on the particle distribution $\rho$
is a convolution operator involving a
matrix-valued kernel $K_h$ defined by
\begin{align*}
K_h(y) &= \nabla_y \left( \Lan \left( \frac{|y|}{h} \right) \; \frac{y}{|y|} \right)
= \frac{1}{h} f_1\left( \frac{|y|}{h} \right) \; I  + \frac{1}{h} f_2\left( \frac{|y|}{h} \right) \; \frac{y}{|y|}  \frac{y^T}{|y|}.
\end{align*}
The coefficients $f_1$ and $f_2$ are given by
\begin{align*}
f_1(z) &= \frac{\Lan \left( z \right)}{z} = \frac{1}{3} - \frac{z^2}{45} + \Ord \left(z^4\right), \\
f_2(z) &= \Lan' \left( z \right) - f_1(z) = - \frac{2 \; z^2}{45} + \Ord \left(z^4\right), \qquad  \text{as} \quad z \to 0.
\end{align*}
$f_1$ and $f_2$ and are analytic functions \cite{marz2016model}. 
Moreover, $K_h$ is a field of symmetric matrices
with pointwise trace $\kappa_h(y) = \trace(K_h(y))$.
The trace kernel can be expressed in terms of the analytic function
\begin{align}\label{eqn:kappa}
f(z) &:= n \; f_1(z) + f_2(z), &
\kappa_h(y) &= \frac{1}{h} f\left(\frac{|y|}{h}\right).
\end{align}
The matrix-valued kernel $K_h$ results from
applying the chain rule in Eq.~\eqref{eqn:signalTrafo}, i.e.
\begin{align}\label{eqn:CoreOp}
s(t) &= \int\limits_{\R^n} \rho(x) K_h(r(t)-x) \; dx \cdot r'(t).
\end{align}
From Eq.~\eqref{eqn:CoreOp} we extract the MPI core operator $A$ 
which acts on $\rho$ via componentwise convolution with $K_h$;
in particular, the MPI core response $A[\rho]$ is given by $A[\rho] = K_h \ast \rho$.
This means that the measured data $s(t)$ is obtained by evaluating 
the MPI core response $A[\rho]$ on the phase space $(r,r')$ of the 
trajectory $r(t)$:
\begin{align}\label{eqn:Split}
s(t) &= A[\rho]( r(t) ) \cdot r'(t) & 
&\text{with}&
A[\rho](y) &= K_h \ast \rho(y).
\end{align}

\section{Methodology and Algorithms}\label{sect:MethodAlgo}
\subsection{Major Stages of the Reconstruction Method}
The proposed reconstruction methods feature two major stages which are 
based on Eq.~\eqref{eqn:Split}.
Under the assumption that a complete sampling of $A[\rho]$ was available, 
i.e., $A[\rho](y)$ were available for all $y$, 
the particle distribution $\rho$ could be reconstructed by deconvolution
from every component of $A[\rho]$ according to the relationship
\begin{align}\label{eqn:ComponentConv}
K_{h,i,j} \ast \rho &= A_{i,j}.
\end{align}
From now on we drop the dependence on $\rho$ and write simply $A$ 
for better readability.
Alternatively, as proposed by the authors in the preliminary 
work \cite{marz2016model},
$\rho$ can be reconstructed by deconvolution relying on
\begin{align}\label{eqn:TraceConv}
\kappa_{h} \ast \rho &= u & 
& \text{with} &
u&:=\trace(A).
\end{align}
A nice property of $\kappa_h$ in contrast to the components $K_{h,i,j}$ 
is that $\kappa_h$ is radially symmetric and positive.
The challenge, on the other hand, is the severe ill-posedness of the 
deconvolution with $\kappa_h$
(cf. \cite[Theorem 3.3]{marz2016model}) because of the kernel's analyticity. 
By exactly the same reasoning, all the deconvolution problems 
involving $K_{h,i,j}$ or $\kappa_h$ are
severely ill-posed, which makes regularization necessary.

The deconvolution problem given by Eq.~\eqref{eqn:TraceConv} is 
actually solved in the second stage, called the \textbf{deconvolution stage}, 
of our method \cite{marz2016model,gapyak2022mdpi,gapyak2023multipatch}.
The reason for this is that, with the usual data acquisition in MPI, we never
get a complete sampling of $A$.
Instead, we obtain a time series of measured data which by Eq.~\eqref{eqn:Split} 
provides by the relationship
\begin{align}\label{egn:MPICorePhase2Signal}
s_l &= A(r_l) v_l, & 
&\text{where} &
s_l &= s(t_l), & 
r_l &= r(t_l), & 
v_l &= r'(t_l),
\end{align} 
the linear actions of samples $A(r_l)$ on velocities $v_l$ along a 
scan trajectory $r(t)$
in terms of the time samples $s_l$.
Thus, the method features a first stage, called the \textbf{MPI core stage}, 
in which the MPI core response $A$ is estimated from 
the time series $s_l$ of measured data. 
The estimate of $A$, resulting from the core stage,
yields then the input data $u := \trace(A)$ for the deconvolution 
problem of Eq.~\eqref{eqn:TraceConv},
which is solved for $\rho$ in the deconvolution stage.

\subsection{MPI Core Stage: Estimation of the MPI Core Response}
The estimation of the core response $A$ is the first stage of our two-stage 
algorithm.
The task here is a data fitting problem employing the relationship 
between the measured data samples and the matrix field $A$ given 
by Eq.~\eqref{egn:MPICorePhase2Signal}.
One special aspect in the MPI-setup is that data fitting has to be done 
on time data and not directly on space data. 
In that regard the task is different from fitting, for example, a spline function.
The signal vectors $s_l = s(t_l)$ are measured at discrete times 
$t_l$, $l \in \{0,\ldots,L-1\}$ 
and constitute the data from which the MPI core response 
$A: \Omega \to \R^{n \times n}$ 
as function of space has to be estimated.
Via the times $t_l$ the position $r_l = r(t_l)$ and velocity $v_l = v(t_l)$ 
vectors correspond to the measured signal vectors $s_l$.
By Eq.~\eqref{egn:MPICorePhase2Signal}, in an ideal scenario, 
$A$ must be such that $s_l - A(r_l) v_l = 0$.
However, in a realistic noisy scenario this exact (over-)fitting 
is not reasonable.
In \cite{gapyak2022mdpi} the authors proposed to estimate $A$ as a 
grid function by minimizing an energy functional $E$ of the form
\begin{align}\label{eqn:EnergyBasic}
E[A] &= \lambda R[A] + F[A].
\end{align}
The grid over the domain $\Omega$ is an $N_x \times N_y$ pixel grid 
and consequently $A$ an $N_x \times N_y \times 2 \times 2$-tensor in a 2D setup. 
In 3D the grid would be  an $N_x \times N_y \times N_z$ voxel grid with
$A$ an $N_x \times N_y \times N_z \times 3 \times 3$-tensor.
The summand $F$ of the energy functional is the data fidelity term 
defined as follows
\begin{align}\label{eqn:fidelity}
F[A] &= \frac{1}{2L} \sum\limits_{l=0}^{L-1} \left| s_l - I[A](r_l) v_l \right|^2.
\end{align}
Here $I[A]$ denotes an interpolation scheme, which is necessary because 
the locations $r_L$ are typically not grid points (i.e. cell centers).
The fidelity term is essentially the sum of the squared residuals w.r.t 
the relationship given in Eq.~\eqref{egn:MPICorePhase2Signal}
using the Euclidean norm $|\cdot|$ on $\R^n$. 
Moreover, the energy in Eq.~\eqref{eqn:EnergyBasic} comprises a regularization term $R$, with regularization 
weight $\lambda > 0$.
The regularizer implicitly selects the search space and,
apart from that, the regularizer controls the gap filling in between data 
locations.

Gap filling is an important issue, because 
typical scanning sequences employed in MPI (cf. \cite{knopp2020openmpidata}) 
use Lissajous trajectories with rather low frequency ratios and sparse sampling
such that the trajectories leave quite big gaps and
many grid cells (especially in higher resolution grids) are never visited 
by the FFP, see Fig.~\ref{fig:liss:scans} (left) for a 2D scanning sequence. 
Here we focus on the 2D setup and aim to reconstruct with higher resolution 
(such as $100 \times 100$ as opposed to current $20 \times 20$ or $40 \times 40$ MPI reconstructions \cite{knopp2020openmpidata})
from scanning data coming from the standard 2D scanning sequence shown 
in Fig.~\ref{fig:liss:scans} (left).
Thus, we have to cope with the fact that the sample locations $r_l$  
give a sparse sampling of the gridded domain $\Omega$ and that
the MPI core response $A$ cannot be determined by data fitting alone. 
Appropriate regularization, which is able to fill the gaps, is needed.

\def\imratio{0.35}
\begin{figure}[t]
\centering
\captionsetup[subfigure]{aboveskip=0pt,belowskip=0pt}
\begin{subfigure}[t]{\imratio\textwidth}
\includegraphics[width=\textwidth]{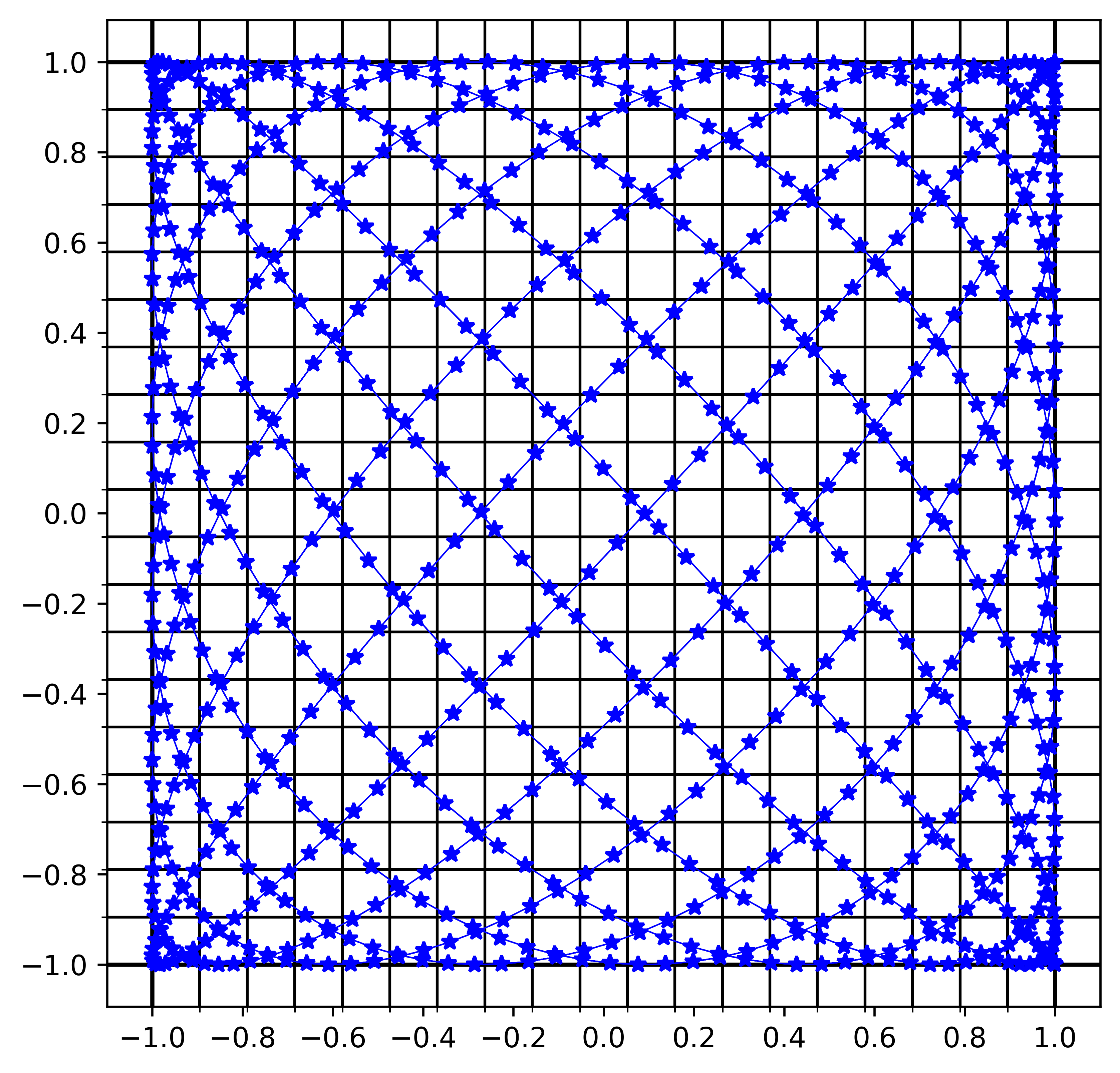}
\end{subfigure}
\hfil
\begin{subfigure}[t]{\imratio\textwidth}
\includegraphics[width=\textwidth]{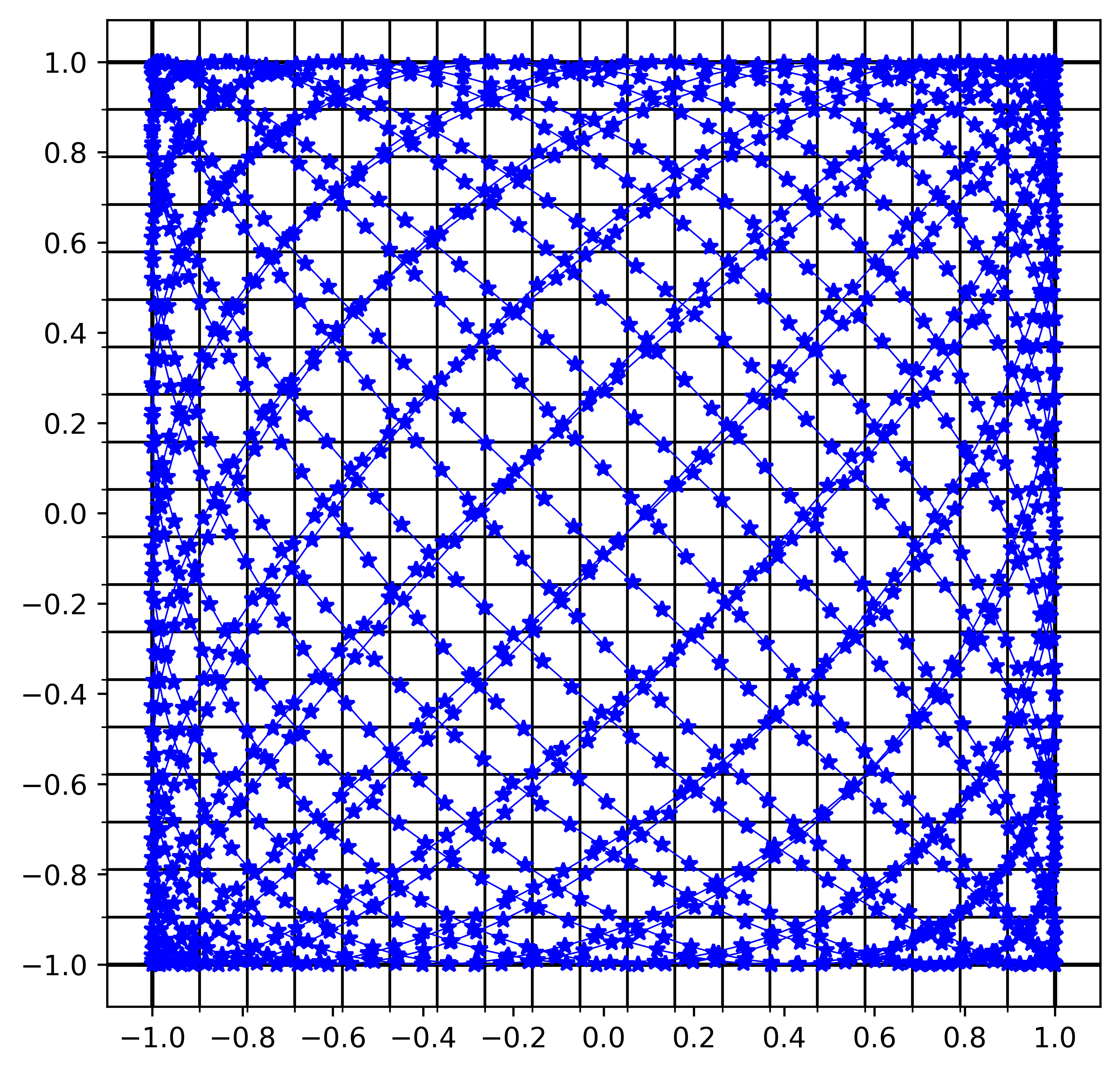}
\end{subfigure}
\caption{\footnotesize 
Scanning sequence: 
Left:
2D Lissajous curve with 16:17 frequency ratio and 1632 samples.
The pixel grid size is $19 \times 19$. 
Both grid and scanning sequence are used for the 2D Open MPI Data \cite{knopp2020openmpidata}.
There are gaps, i.e. pixels never visited by the
FFP which moves on the Lissajous curve.
When working with finer grids, the gap effect becomes more pronounced. 
Right: a denser scan
obtained by merging the scan left with 
same rotated by $90^\circ$ as suggested by \cite{gapyak2022mdpi}.
}
\label{fig:liss:scans}
\end{figure}

\subsubsection{MPI Core Stage: 1st Order Regularization Using 
Eigenfunctions}\label{sect:PrelimAlgo}

In the preliminary work \cite{gapyak2022mdpi,gapyak2023multipatch} 
we employed the regularizer
\begin{align}\label{eqn:Regularizer1}
R_1[A] &= \frac{1}{2\; N_x N_y} \; \|D A\|_{2}^2 
\end{align}
where $D$ denotes the approximation of the gradient by finite differences.
For the interpolation operation $I[A]$ in the data fidelity term we used 
bi-quadratic interpolation. 

Here, we take a different route: we will approximate the MPI core response
$A:\Omega \to \R^{2\times 2}$ with analytic functions.
This is because, as discussed in Sect.~\ref{sect:CoreOp}, $A$ is 
supposed to be analytic.
In particular, we use the cosine expansion     
\begin{align}\label{eqn:expansion}
A(x,y) &= \sum\limits_{m \in \N_0^2} \hat{A}_m u_m(x,y), &
& (x,y) \in \Omega = [-1,1]^2,
\end{align}
with matrix-valued coefficients $\hat{A}_m \in \R^{2\times 2}$. 
The cosine basis is given by
\begin{align}\label{eqn:CosBasis}
u_m(x,y) &= c_m \cos \left( \frac{\pi m_1}{2} (x+1)\right)
 \cos \left( \frac{\pi m_2}{2} (y+1)\right), 
\qquad m = (m_1,m_2) \in \N_0^2,  
\end{align}
where the constants $c_m$ are such that the $L_2$-Norm of $u_m$ equals $1$.

The reason for choosing the cosine basis is that the regularizer 
of Eq.~\eqref{eqn:Regularizer1} 
is a discrete version of
\begin{align}\label{eqn:RegularizerCont1}
R_1[A] &= \sum\limits_{p=1}^2 \sum\limits_{q=1}^2 \frac{1}{2 |\Omega|} \int\limits_\Omega \left| \nabla A^{p,q}(x,y) \right|^2 \; d(x,y).
\end{align}
The underlying operator of $R_1$, by virtue of the calculus of variations, 
is the negative Laplacian with Neumann-Zero boundary conditions.
The corresponding eigenbasis, i.e., the solutions of
\begin{align}\label{eqn:EigenLapNat}
-\laplace u &= \mu u \qquad \text{in} \; \Omega, & 
\partial_{\nu} u &= 0 \qquad \text{on} \; \partial \Omega,
\end{align}
are the cosine functions stated in Eq.~\eqref{eqn:CosBasis} with eigenvalues
\begin{align}\label{eqn:EigvalsLaplace}
\mu_m &
= \frac{\pi^2}{4} \left(m_1^2 + m_2^2\right)\;.
\end{align}
The Neumann-zero boundary conditions are natural w.r.t. 
the calculus of variations, that means that the energy functional
\begin{align}
E_1[A] &= \lambda \sum\limits_{p=1}^2 \sum\limits_{q=1}^2 \frac{1}{2 |\Omega|} \int\limits_\Omega \left| \nabla A^{p,q}(x,y) \right|^2 \; d(x,y)
+ \frac{1}{2L} \sum\limits_{l=0}^{L-1} \left| s_l - A(r_l) v_l \right|^2
\end{align} 
is well-defined on the space
\begin{align}
X_1 &:= \left\{ A:\Omega \to \R^{2 \times 2} \;:\; A^{p,q} \in H^1(\Omega) \;,\; \left. \partial_{\nu} A^{p,q} \right|_{\partial \Omega} = 0\right\}.
\end{align}
$H^1(\Omega)$ denotes the Sobolev space of $L^2$-functions with first order weak derivatives in the set $L^2(\Omega)$ of square-integrable functions.

Having $A$ expanded in terms of the cosine functions, we can rewrite the 
energy $E$ in terms of the
sequence $\hat{A} = (\hat{A}_m)_{m\in \N_0^2}$ of expansion coefficients.
Plugging the expansion given in Eq.~\eqref{eqn:expansion} into the 
regularizer $R_1$, we obtain
\begin{align}\label{eqn:RegularizerR1}
R_1[A] &= \sum\limits_{p=1}^2 \sum\limits_{q=1}^2 \frac{1}{2 |\Omega|} 
\sum\limits_{m \in \N_0^2} \mu_m \;\left(\hat{A}^{p,q}_m\right)^2 
= \frac{1}{2 |\Omega|}  \sum\limits_{m \in \N_0^2} \mu_m \; \| \hat{A}_m \|_F^2   
=: \hat{R}_1[\hat{A}]
\end{align}
where $\|\hat{A}_m\|_F$ is the Frobenius norm of the 
$2 \times 2$ matrix $\hat{A}_m$.
Thus the energy is given by
\begin{align}\label{eqn:FunctionalE1}
E_1[A] &= \frac{\lambda}{2 |\Omega|}  \sum\limits_{m \in \N_0^2} \mu_m \; \| \hat{A}_m \|_F^2 
+ \frac{1}{2L} \sum\limits_{l=0}^{L-1} \left| s_l - \sum\limits_{m \in \N_0^2} \hat{A}_m v_l \; u_m(r_l)  \right|^2
=: \hat{E}_1[\hat{A}]
\end{align}

For the purpose of discretization we truncate the expansion to employ 
only the first $N \times M$  basis functions
\begin{align}\label{eqn:expansionFin}
A(x,y) &= \sum\limits_{k=0}^{N-1} \sum\limits_{l=0}^{M-1} \hat{A}_{k,l} \; u_{k,l}(x,y).
\end{align}
This function $A$ is related to a grid function defined at cell center 
points via
\begin{align}
A_{i,j} &:= A(x_i,x_j) = \sum\limits_{k=0}^{N-1} \sum\limits_{l=0}^{M-1} \hat{A}_{k,l} \; u_{k,l}(x_i,x_j).
\end{align} 
That means that the finite sequences $A_{i,j}$ and $\hat{A}_{k,l}$ 
are related via the discrete cosine transform
when $N=N_x,M=N_y$ are tied to the grid size $N_x,N_y$.

After truncation we find the minimizer of the quadratic functional 
$\hat{E}_1[\hat{A}]$ with $\hat{A} \in \R^{N \times M \times 2 \times 2}$
by solving the gradient system $\nabla \hat{E}_1[\hat{A}] = 0$, 
which is a system of linear equations
with a symmetric positive definite system matrix given by the 
Hessian of $\hat{E}_1$, with 
the conjugate gradients (CG) method.  

\subsubsection{MPI Core Stage: Novel Algorithm with 2nd Order Regularization}\label{sect:ImprovedAlgo}

The preliminary algorithm has limitations, which are demonstrated 
in Sect.~\ref{sect:Exp1} with experiment 1.
To improve the performance of the overall reconstruction of the 
particle density $\rho$, we propose
to exploit second order regularization in the MPI core stage.
That is, we employ the regularizer
\begin{align}\label{eqn:RegularizerR2}
R_2[A] &= \sum\limits_{p=1}^2 \sum\limits_{q=1}^2 \frac{1}{2 |\Omega|} \int\limits_\Omega \left| \laplace A^{p,q}(x,y) \right|^2 \; d(x,y),
\end{align}
and the energy functional becomes
\begin{align}\label{eqn:FunctionalE2} 
E_2[A] &= \lambda \sum\limits_{p=1}^2 \sum\limits_{q=1}^2 \frac{1}{2 |\Omega|} \int\limits_\Omega \left| \laplace A^{p,q}(x,y) \right|^2 \; d(x,y)
+ \frac{1}{2L} \sum\limits_{l=0}^{L-1} \left| s_l - A(r_l) v_l \right|^2.
\end{align} 
This energy functional is well-defined on the space
\begin{align}
X_2 &:= \left\{ A:\Omega \to \R^{2 \times 2} \;:\; A^{p,q} \in H^2(\Omega) \;,\; \left. \partial_{\nu} A^{p,q} \right|_{\partial \Omega} = 0\right\},
\end{align}
where $H^2(\Omega)$ denotes the Sobolev space of $L^2$-functions with first and second order weak derivatives in the set $L^2(\Omega)$ of square-integrable functions.
By keeping the cosine expansions of Eq.~\eqref{eqn:expansion}, which are
the eigenfunctions of the Laplacian,
we obtain $R_2$ also in terms of the expansion coefficients, i.e.,
\begin{align}\label{eqn:Regularizer2}
R_2[A] &= \frac{1}{2 |\Omega|}  \sum\limits_{m \in \N_0^2} \mu_m^2 \; \| \hat{A}_m \|_F^2   
=: \hat{R}_2[\hat{A}].
\end{align}
Accordingly, we obtain the energy $\hat{E}_2[\hat{A}]$ expressed 
in terms of the expansion coefficients.
Note, that $\hat{R}_1[\hat{A}]$ and $\hat{R}_2[\hat{A}]$ are 
both essentially weighted $l^2$-norms,
but differ greatly in the weights. The weights of $\hat{R}_2[\hat{A}]$
are the squared eigenvalues $\mu_m$ of the Laplacian from 
Eq.~\eqref{eqn:EigvalsLaplace}
while $\hat{R}_1[\hat{A}]$ has directly the eigenvalues as weights.
Hence, higher modes corresponding to higher frequencies are penalizes 
much stronger
when employing $\hat{R}_2[\hat{A}]$ as regularizer.

\subsection{Deconvolution Stage: Reconstruction of the Particle Distribution}
In the deconvolution stage, the second stage of the two-stage method, 
we reconstruct the particle density $\rho$.
Given the estimation of $A$ from the MPI core stage, 
we solve the ill-posed deconvolution problem by energy minimization.
The basic form of the energy functional to be minimized is given by
\begin{align}\label{eqn:FunctionalJtraceGeneric}
J[\rho] &= \mu \; R_D[\rho] + \| \kappa_{h} \ast \rho - u \|_{L^2(\Omega)}^2
\end{align}
where $\mu >0$ is the regularization weight, $R_D$ the regularization prior, 
$u = \trace(A)$ is the data,
and $\kappa_{h}$ the kernel defined in Eq.~\eqref{eqn:kappa}.
The fidelity term here is the residual of the convolution relation~\eqref{eqn:TraceConv}
which $\rho$ should satisfy w.r.t. the trace of $A$.

In our earlier work \cite{marz2016model,gapyak2022mdpi,gapyak2023multipatch,gapyak2025fast} we have 
suggested and studied different regularization priors.
Here, we use our most recent approach from \cite{gapyak2025fast},
which is a Plug-and-Play prior \cite{venkatakrishnan2013pnp}.
We summarize the approach briefly.
Half Quadratic Splitting (HQS) considers the constrained minimization problem
\begin{align}\label{eq:dec:min:2}
\rho &= \arg\min_{\rho_1 \, , \rho_2}\lVert\kappa_h * \rho_1 - u\rVert_2^2 + \mu R_D [\rho_2]\quad\text{s.t. } \rho_1 = \rho_2\, ,
\end{align}
which is equivalent to the free minimization of $J$,
and minimizes the Lagrangian w.r.t. the two variables $\rho_1$ and $\rho_2$
in an alternating fashion:
\begin{align}
	\rho_1^{k+1} & = \arg\min_{\rho_1} \lVert \kappa_h * \rho_1 -u\rVert_2^2 +  \nu_k\left\lVert \rho_1-\rho_2^k \right\rVert_2^2 \label{eq:dec:sub:tik0}\\
	\rho_2^{k+1} & = \arg\min_{\rho_2} \mu R_D [\hat{\rho}_2]+\nu_k\left\lVert \rho_1^{k+1}-\rho_2 \right\rVert_2^2 . \label{eq:dec:sub:den:gauss}
\end{align}
Eq.~\eqref{eq:dec:sub:den:gauss} is a denoising problem in variational form,
where Gaussian noise of noise level $\sqrt{\mu /\nu_k}$ is to be removed from
the variable $\rho_1^{k+1}$. 
Following the Plug-and-Play ansatz \cite{venkatakrishnan2013pnp}, 
Eq.~\eqref{eq:dec:sub:den:gauss} is replaced with a denoiser. 
We use the denoiser from \cite{Zhang2022pnp}, a pre-trained UNet-based
Gaussian denoiser called DRUNet.
The DRUNet can take into account general noise level maps.
Here, the noise level map is the constant $\sqrt{\mu /\nu_k}$.
The resulting algorithm, based on Eq.~\eqref{eq:dec:sub:tik} and \eqref{eq:dec:sub:den}, alternates between a simpler Tikhonov-regularized
reconstruction step involving the data-fidelity and a Gaussian denoising step
\begin{align}
\rho_1^{k+1} & = \arg\min_{\rho_1}\lVert u-C_h \rho_1\rVert_2^2 +  \nu_k\left\lVert \rho_1-\rho_2^k \right\rVert_2^2 \label{eq:dec:sub:tik}\\
\rho_2^{k+1} & = \mathrm{Denoiser}\left (\rho_1^{k+1}\, , \sqrt{\mu /\nu_k}\right ) . \label{eq:dec:sub:den}
\end{align}
$C_h$ denotes the discretized version of the convolution operator with kernel $\kappa_h$.
As discussed in \cite{gapyak2025ell} the noise level $\sigma_{k+1}=\sqrt{\mu /\nu_k}$ 
is inversely related to the Tikhonov parameter $\nu_k$. 
Therefore, we have only $\nu_0$ as a (hyper)parameter, because $\sigma_{k+1}$ 
is estimated from the iterate $\rho_1^{k+1}$ with the root of the variance 
\begin{align}
	\sigma_{k+1} &= \sqrt{\mathrm{Var}(\rho_1^{k+1})}\;.
\end{align}
This estimate of $\sigma_{k+1}$ has proven to be reliable and yields good reconstructions \cite{gapyak2025ell}.

The discrete deconvolution problem of Eq.~\eqref{eq:dec:sub:tik} is formally a regularized inversion of the Toepliz matrix $C_h$. 
But, because $C_h$ is not sparse we make use of the FFT and apply direct multiplication in Fourier domain instead of matrix multiplication by $C_h$. 
The minimization problem of Eq.~\eqref{eq:dec:sub:tik}  is solved with
the CG method on the Euler-Lagrange equations.

The focus of this paper is the MPI core stage, that is for the deconvolution stage
we employ the method of \cite{gapyak2025fast} outlined above.
In Sect.~\ref{sect:Exp1} and Sect.~\ref{sect:Exp2}, we illustrate with 
experiments that the proposed change of the MPI core stage is indeed helpful
and yields enhanced overall reconstruction quality.

\section{Theoretical Foundations of the Approach in the MPI Core Stage}\label{sect:Theory}
The motivation for the choice of regularizers in Sect.~\ref{sect:MethodAlgo} 
is smoothness of the core response $A$, 
since the data for the deconvolution stage of the 
reconstruction scheme is $u = \trace(A)$ the pointwise trace of $A$.
In an ideal scenario $u$ and $A$ should be smooth because they satisfy 
$\kappa_h \ast \rho = u$ and $K_h \ast \rho = A$ 
with analytic convolution kernels $\kappa_h$ and $K_h$, respectively. 
Thus a smoothness supporting regularizer such as the $L^2$-norm of the gradient,
defined in Eq.~\eqref{eqn:RegularizerCont1}, seems reasonable and is our starting choice 
to be used with $E$.
However, it turns out that when using $R_1$ with sparse data, 
then $u$ looks very spiky (see Fig.~\ref{subfig:trace1sparse}).
We note that the regularizer $R_1$ relates also to the energy of a membrane.
Inspired by smoothing in geometry processing \cite{stein2018natural}
more reasonable choices for regularizers 
to enforce higher smoothness (see Fig.~\ref{subfig:trace2sparse})
are the $L^2$-norm of the Laplacian, defined in Eq.~\eqref{eqn:Regularizer2},
and the $L^2$-norm of the Hessian, given by
\begin{align}\label{eqn:Regularizer3}
R_3[A] &= \sum\limits_{p=1}^n \sum\limits_{q=p}^n \frac{1}{2 |\Omega|} \int\limits_\Omega \left| \hessian A^{p,q}(x) \right|_F^2 \; dx.
\end{align}
Both regularizers $R_2$ (Eq.~\eqref{eqn:RegularizerR2}) and $R_3$ have the biharmonic operator 
(or Bi-Laplacian) as underlying linear operator, via the calculus of variations, but with different
options to impose enforced or natural boundary conditions.
Moreover, both regularizers relate to the energy of a thin plate.

Knowing the underlying linear operator (with boundary conditions),
which describes the gradient of the regularizer,
it makes sense to expand the solution of the minimization problem
w.r.t. the eigenfunctions. 
If the eigenfunctions can be computed analytically, they are very helpful to
derive a simpler description of the regularizer.
However, the concrete eigenfunctions depend on the domain $\Omega$ and the chosen boundary conditions. 

In our application the domain $\Omega = [-1,1]^2$ is an axis-parallel square
where tensor-product eigenfunctions are in principle possible.
For the negative Laplacian with (natural) Neumann-zero boundary conditions,
the tensor-product cosine functions from Eq.~\eqref{eqn:CosBasis} are exactly these eigenfunctions.
In image processing the (natural) Neumann-zero boundary conditions are the default, 
when working with $R_1$ (Eq.~\eqref{eqn:RegularizerCont1}).
But regarding a membrane attached to a frame, Dirchlet-zero boundary conditions are reasonable.
For the negative Laplacian with (enforced) Dirchlet-zero boundary conditions,
the tensor-product sine functions from Eq.~\eqref{eqn:SinBasis} are the eigenfunctions.  

For the regularizers $R_2$ (Eq.~\eqref{eqn:RegularizerR2}) and $R_3$ (Eq.~\eqref{eqn:Regularizer3}) and the underlying Bi-Laplacian the situation is more complex.
Because of the higher order there are two sets of boundary conditions to be chosen.
It turns out that tensor-product eigenfunctions are not necessarily available for all sets of boundary conditions.
Since natural boundary conditions worked so well with $R_1$, it may appear reasonable 
to have all natural boundary conditions also with $R_2$.
Interestingly, we will see that in this case only for the eigenvalue $\mu=0$ tensor-product eigenfunctions are 
available (Thm.~\ref{theo:EigBiLapAllNatSep}), i.e., they are all in the kernel. 
Furthermore, in this case the underlying operator has infinite dimensional kernel
containing all harmonic functions
(Thm.~\ref{theo:BiLapAllNatKernel}) and thus the regularizer is not helpful.
Therefore, it makes sense to enforce at least partially certain boundary conditions.
Enforcing Neumann-zero boundary conditions and leaving the remaining boundary conditions to come in naturally 
(in the sense of the calculus of variations) yields again the tensor-product cosine functions.
We discuss also, for completeness, the case of enforcing Dirichlet-zero boundary conditions
while leaving remaining boundary conditions to come in naturally
which yields once more the tensor-product sine functions.
This latter case can be encountered in mechanics when considering a simply supported thin plate
and constructing the Navier solution.
We consider also all enforced boundary conditions, as these are a standard case 
when considering a clamped thin plate. In this case, there a no tensor-product eigenfunctions 
(Thm.~\ref{theo:BiLapEigenAllEnf}).
 
In geometry processing Hessian regularization with $R_3$ is also common.
Even though the regularizers $R_2$ and $R_3$ look different in terms of their definition
and $R_3$ appears to be more complex, we reveal  
that in our setup ($\Omega$ a square, enforced Neumann-zero boundary conditions) 
$R_2$ and $R_3$ are actually identical (Thm.~\ref{theo:LapReguEquivHessRegu}
and Thm.~\ref{theo:R2R3sameEig}).
Therefore, we considered only $R_2$ in Sect.~\ref{sect:ImprovedAlgo}.
In the following we discuss systematically 
the eigenvalue problems of the Laplacian and Bi-Laplacian
w.r.t. the existence of tensor-product eigenfunctions.

\subsection{Eigenfunctions of the Laplacian}
We collect the results for the Laplacian, as they will reappear when we consider the Bi-Laplacian.
For a scalar-valued function $u:\Omega \to \R$ the regularizer $R_1$ is given by
\begin{align*}
R_1[u] &= \frac{1}{2 |\Omega|} \int\limits_\Omega \left| \nabla u(x) \right|^2 \; dx
\end{align*}
with variational derivative $\delta R_1[u] \; \varphi = \left. \frac{d}{dh} R_1[u+h\varphi] \right|_{h=0}$
given by
\begin{align}\label{eqn:LapBI}
\delta R_1[u] \; \varphi &= 
\frac{1}{|\Omega|} \left( \int\limits_\Omega -\laplace u(x) \varphi(x) \; dx + \int\limits_{\partial \Omega} \partial_{\nu} u(x) \; \varphi(x) \; dS(x) \right). 
\end{align}
Here $\nu$ denotes the exterior normal to the boundary $\partial\Omega$ of $\Omega$
and $\partial_{\nu} u$ the derivative in direction $\nu$.

\textbf{Laplacian with natural boundary conditions.}
If there are no enforced boundary conditions on $u$, then $\varphi$ is also free on the boundary.
That means for the corresponding energy minimization that $\partial_{\nu} u$ must vanish
for $\delta E[u] \; \varphi$ to vanish.
This results in the natural boundary conditions (which are often the default 
in image and geometry processing) to complete the eigenvalue problem,
which is stated in Eq.~\eqref{eqn:EigenLapNat}.
Because $\Omega$ is a square in 2D the eigenfunctions can be solved for 
by separation of variables. We obtain exactly the cosine basis of Eq.~\eqref{eqn:CosBasis};
cf.~\cite{courant2004methods}.
Moreover, we see that for the task of minimization the space
\begin{align}\label{eqn:SpaceX1Scalar}
X_1 &:= \left\{ u:\Omega \to \R \mid \Omega = [-1,1]^2, \; u \in H^1(\Omega), \; \partial_\nu u |_{\partial \Omega} = 0\right\},
\end{align}
is the reasonable domain for $R_1$ when there are no additional conditions.
These eigenfunctions are also orthonormal w.r.t. the scalar product on the space $L^2(\Omega)$.
If $u \in X_1$ is expanded as 
\begin{align}\label{eqn:CosExpansion}
u(x,y) &= \sum\limits_{m \in \N_0^2} \hat{u}_m \; u_m(x,y)
\end{align}
using the cosine basis $u_m$ of Eq.~\eqref{eqn:CosBasis},
then we have
\begin{align}
R_1[u] &= \frac{1}{2 |\Omega|} \int\limits_\Omega \left| \nabla u(x) \right|^2 \; dx =
\frac{1}{2 |\Omega|} \sum\limits_{m \in \N_0^2} \mu_m \;\left(\hat{u}_m\right)^2 
\end{align}
with $|\Omega| = 4$.
Hence, $R_1$ is a diagonal quadratic form on the sequence $(\hat{u}_m)_{m \in \N_0^2}$ of coefficients.
In the case of matrix-valued functions with components $A^{p,q} \in X_1$ 
we expand the components accordingly and obtain Eq.~\eqref{eqn:RegularizerR1}

\textbf{Laplacian with enforced boundary conditions.}
For the sake of completeness, we mention also the use of $R_1$ with enforced Dirichlet-zero boundary conditions,
i.e. $u=0$ on $\partial \Omega$.
Such enforced boundary conditions are considered in mechanics, e.g., in the case of a membrane held fixed at the boundary.
The side effect of the enforced boundary is that $\varphi$ must also vanish on $\partial \Omega$
for $u+h\varphi$ to be admissible. Thus, the boundary integral in Eq.~\eqref{eqn:LapBI} vanishes automatically.
The eigenvalue problem with Dirichlet-zero boundary conditions now reads as
$-\laplace v = \mu v$ in $\Omega$ and $v=0$ on $\partial \Omega$.
And in this setup the eigenfunctions on the square $\Omega$
are given by (cf. \cite{courant2004methods, berard2015dirichlet}) 
\begin{align}\label{eqn:SinBasis}
v_m(x,y) &= c_m \sin \left( \frac{\pi m_1}{2} (x+1)\right)
 \sin \left( \frac{\pi m_2}{2} (y+1)\right), 
\qquad m = (m_1,m_2) \in \N^2.
\end{align} 
The eigenvalues are still those of Eq.~\eqref{eqn:EigvalsLaplace}, but note that this time
the indices $m_1,m_2$ are non-zero.

\subsection{Eigenfunctions of the Bi-Laplacian}
When we work with the regularizer $R_2$ there are more options to impose boundary conditions
because of the second order derivative under the integral.
For scalar-valued $u:\Omega \to \R$, we have
\begin{align*}
R_2[u] &= \frac{1}{2|\Omega|} \int\limits_\Omega \left(\laplace u \right)^2 \; dx
\end{align*} 
with variational derivative
\begin{align}
\delta R_2[u] \; \varphi &= \left. \frac{d}{dh} R_2[u+h\varphi] \right|_{h=0} =
\left. \frac{1}{2|\Omega|} \frac{d}{dh} \int\limits_\Omega \left(\laplace u + h \laplace \varphi \right)^2 \; dx \right|_{h=0}
\nonumber \\
&=\frac{1}{|\Omega|} \int\limits_\Omega \laplace u  \laplace \varphi \; dx
= \frac{1}{|\Omega|} \left( \int\limits_\Omega \laplace^2 u \; \varphi \;dx 
+ \int\limits_{\partial \Omega} \laplace u \; \partial_{\nu} \varphi - \partial_{\nu} \laplace u \;  \varphi \; dS(x) \right). \label{eqn:BiLapBI}
\end{align}
In the following part of this section, we are going to discuss the eigenvalue problems with different 
choices of boundary conditions in more detail. Interestingly, the reference \cite{courant2004methods},
often quoted w.r.t. to this type of problems,
does not cover the cases of other than all enforced boundary conditions,
i.e., only the case of the clamped thin plate.

\textbf{Bi-Laplacian with all natural boundary conditions.}
Because with the regularizer $R_1$ (and the Laplacian) 
no boundary conditions are explicitly enforced and
we let them enter the scene naturally, it seems reasonable, at first glance,
to do the same when using the regularizer $R_2$.
Without enforced boundary conditions, $\varphi$ as well as $\partial_{\nu} \varphi$ are free, and 
the all natural boundary conditions are $\laplace u=0$ and $\partial_{\nu} \laplace u=0$
to make the boundary integral of Eq.~\eqref{eqn:BiLapBI} zero.
With that the eigenvalue problem reads as
\begin{align}\label{eqn:EigBiLapAllNat}
\laplace^2 u &= \mu u \qquad \text{in} \; \Omega, \qquad & 
\laplace u &= 0 \quad \text{and} \quad
\partial_{\nu} \laplace u = 0 \quad \text{on} \; \partial \Omega.
\end{align}

\begin{theorem}\label{theo:EigBiLapAllNatSep}
The eigenvalue problem of Eq.~\eqref{eqn:EigBiLapAllNat} has
separable solutions only for $\mu = 0$.
In the case of $\mu = 0$ there is a family of separable solutions,
that is non-linear in one of the parameters.
\end{theorem}

\begin{proof}
Aiming at separable solutions on the square $\Omega$,
the ansatz $u(x,y) = \alpha(x) \; \beta(y)$ implies
$\laplace^2 u = \partial_x^4 u + 2 \partial_x^2 \partial_y^2 u + \partial_y^4 u =  \alpha'''' \beta + 2 \alpha'' \beta'' + \alpha \beta''''$
and the problem becomes
\begin{align}
 \alpha'''' \beta + 2 \alpha'' \beta'' + \alpha \beta'''' &= \mu \alpha \beta.
\end{align}
Since $u$ is not the zero function, we convert it to
\begin{align}\label{eqn:EigProbBiLapSep}
 \frac{\alpha''''}{\alpha}  + 2 \frac{\alpha''}{\alpha} \frac{\beta''}{\beta} + \frac{\beta''''}{\beta} &= \mu.
\end{align}
Because $\alpha$ depends only on $x$ and $\beta$ only on $y$ while $\mu$ is a constant, we obtain
\begin{align*}
\left( \frac{\alpha''}{\alpha}\right)' \left(\frac{\beta''}{\beta}\right)'  &= 0
\end{align*}
after differentiating Eq.~\eqref{eqn:EigProbBiLapSep} w.r.t. $x$ and then w.r.t. $y$.
So, at least one of the two factors must be the zero-function, i.e., we have
\begin{align}\label{eqn:Alpha}
\left( \frac{\alpha''}{\alpha}\right)' &= 0 &
& \implies &
\alpha'' &= c \alpha
\end{align}
with a constant $c$. By employing the boundary conditions on $u$, we obtain further conditions on $\alpha$.
By the boundary condition $\laplace u=0$, we have at $x = -1$
\begin{align*}
\laplace u(-1,y) &= \alpha''(-1) \; \beta(y) + \alpha(-1) \; \beta''(y) =  \alpha(-1) \; (\beta''(y) + c \beta(y))
= 0 \qquad \text{for all} \; y.
\end{align*}
By the second boundary condition $\partial_{\nu} \laplace u=0$, we have in addition
\begin{align*}
\partial_{\nu} \laplace u(-1,y) &= -\partial_x \laplace u(-1,y) = -\alpha'(-1) \; (\beta''(y) + c \beta(y))
= 0 \qquad \text{for all} \; y.
\end{align*}
If we assume for now, that $\beta$ does not satisfy the ordinary differential equation $\beta'' = -c \beta$,
then we must have $\alpha(-1) = 0$ and $\alpha'(-1) = 0$.
In this case $\alpha$ must solve the initial value problem
\begin{align}\label{eqn:IVP}
\alpha'' &= c \alpha \;, &
\text{with} \qquad \alpha(-1) &= 0 \;, &
\text{and} \qquad \alpha'(-1) &= 0.
\end{align}
The unique solution to this problem is $\alpha \equiv 0$, leading to $u \equiv 0$, 
which is not admissible as eigenfunction.
The only possibility is that we have also $\beta'' = -c \beta$.
But then the Laplacian
\begin{align*}
\laplace u &= \alpha'' \beta + \alpha \beta'' = c \alpha \beta + \alpha \beta''
= \alpha (\beta'' + c \beta ) = 0
\end{align*}
vanishes identically and thus also the Bi-Laplacian. 
Hence, a separable solution is possible only for the eigenvalue $\mu = 0$.
For eigenvalues $\mu > 0$, there are no separable solutions.

Moreover, for the eigenvalue $\mu = 0$ we have infinitely many separable solutions. 
If we choose $c = -\omega^2$ negative in Eq.~\eqref{eqn:Alpha}, then
\begin{align}\label{eqn:KernelSep}
u(x,y) =& \; (a_1 \cos(\omega x) + a_2\sin( \omega x)) \; 
(b_1\cosh( \omega y) + b_2 \sinh(\omega y)) \nonumber \\
=& \; c_1 \cos(\omega x) \cosh( \omega y)
+ c_2 \sin( \omega x) \cosh( \omega y) \\
& \quad + c_3 \cos(\omega x) \sinh( \omega y)
+ c_4 \sin( \omega x) \sinh( \omega y)\nonumber
\end{align}
is a solution for any choice of $\omega$. This family of separable solutions
is non-linear w.r.t. the parameter $\omega$.
If we choose $c = \omega^2$ positive in Eq.~\eqref{eqn:Alpha},
then
\begin{align*}
u(x,y) &= (a_1 \cosh(\omega x) + a_2\sinh( \omega x)) \; 
(b_1\cos( \omega y) + b_2 \sin(\omega y))
\end{align*}
is also a solution for any choice of $\omega$.
If we choose $c=0$ in in Eq.~\eqref{eqn:Alpha}, then
\begin{align*}
u(x,y) &= (a_1  + a_2 x) \; (b_1 + b_2 y)\;.
\end{align*}
\end{proof}

This result implies that the Bi-Laplacian with all natural boundary conditions has a huge kernel.
More precisely, we have the following theorem.

\begin{theorem}\label{theo:BiLapAllNatKernel}
The kernel of the Bi-Laplacian with all natural boundary conditions,
i.e., the set of solutions of Eq.~\eqref{eqn:EigBiLapAllNat} for $\mu=0$,
consists of all harmonic functions on $\Omega$.
The kernel contains a countable orthogonal system and is thus infinite-dimensional.
\end{theorem}

\begin{proof}
The set of all harmonic functions is given by the solutions of $\laplace u = 0$ on $\Omega$
(without any boundary conditions).
Any harmonic function satisfies certainly Eq.~\eqref{eqn:EigBiLapAllNat} for $\mu=0$.
The other way round, if $u$ is a solution of Eq.~\eqref{eqn:EigBiLapAllNat} for $\mu=0$,
we employ integration by part as in Eq.~\eqref{eqn:BiLapBI}
with $\varphi=u$, and obtain 
$\|\laplace u\|_{L^2(\Omega)}^2 = \int\limits_\Omega (\laplace u)^2 \; dx = 0.$
Thus, we must have $\laplace u = 0$ on $\Omega$. This proves the first claim.
For the second claim, we pick from the separable solutions of Eq.~\eqref{eqn:KernelSep}
the following subset
\begin{align*}
u_n(x,y) &= \cos\left( \frac{\pi n}{2} (x+1)\right) \cosh\left( \frac{\pi n}{2} (y+1)\right)\;, \qquad n \in \N_0.
\end{align*}
This is an orthogonal system, which satisfies $\laplace u_n = 0$.
\end{proof}

\textbf{Bi-Laplacian with partially enforced boundary conditions.}
First, we discuss the case we are interested in:
we enforce only the boundary condition $\partial_{\nu} u=0$, 
because this is the typical choice in image and geometry processing (see also \cite{stein2018natural}).
Note that in the first order scenario this boundary condition appears naturally,
but in the current context it is enforced.
With this we must also have $\partial_{\nu} \varphi=0$ for $u+h\varphi$ to be admissible.
Hence, the first summand of the boundary integral in Eq.~\eqref{eqn:BiLapBI} is zero.
The remaining second boundary condition comes naturally (in the spirit of the calculus of variations),
i.e., we must have $\partial_{\nu} \laplace u = 0$ for the second summand of the 
boundary integral in Eq.~\eqref{eqn:BiLapBI} to vanish.
Now, the corresponding eigenvalue problem is given by
\begin{align}\label{eqn:EigBiLapR2}
\laplace^2 u &= \mu u \qquad \text{in} \; \Omega, \qquad & 
\partial_{\nu} u &= 0 \quad \text{and} \quad
\partial_{\nu} \laplace u = 0 \quad \text{on} \; \partial \Omega.
\end{align}
This problem admits separable solutions on the square $\Omega$.
The eigenbasis is exactly the cosine basis $u_m$ from Eq.~\eqref{eqn:CosBasis} with associated eigenvalues $\mu_m^2$,
which can be verified directly.

Employing the expansion of $u \in H^2(\Omega)$ from Eq.~\eqref{eqn:CosExpansion}, we have
\begin{align}
-\laplace u(x) &= \sum\limits_{m \in \N_0^2} \hat{u}_m \; (-\laplace u_m(x))
= \sum\limits_{m \in \N_0^2} \hat{u}_m \; \mu_m \; u_m(x).
\end{align}
Again by orthonormality of the eigenfunctions we have
\begin{align}\label{eqn:Reg2CosBasis}
\int\limits_\Omega \left( \laplace u(x) \right)^2 \; dx &=
\sum\limits_{m \in \N_0^2} \mu_m^2 \;\left(\hat{u}_m\right)^2.
\end{align}
Hence, $R_2$ is a diagonal quadratic form on the sequence $(\hat{u}_m)_{m \in \N_0^2}$ of coefficients.
In the case of matrix-valued functions $A$ this allows us to express the regularizer $R_2[A]$ as in Eq.~\eqref{eqn:Regularizer2}.

The other option which comes to mind for working with partially enforced boundary conditions is
to enforce again Dirichlet-zero boundary conditions, i.e. $u=0$.
As before, this implies $\varphi=0$ on the boundary and thus the second summand of the boundary integral in Eq.~\eqref{eqn:BiLapBI} vanishes automatically.
The remaining natural boundary condition is thus $\laplace u = 0$.
This particular set of boundary conditions comes up, e.g., with the problem of a simply supported thin plate in mechanics.
Here, the eigenproblem reads as
$\laplace^2 u = \mu u$ in $\Omega$ 
with $u = 0$ and $\laplace u = 0$ on $\partial \Omega$.
This problem admits also separable solutions on the square $\Omega = [-1,1]^2$
which are exactly the sine basis $v_m$ from Eq.~\eqref{eqn:SinBasis} 
with eigenvalues $\mu_m^2$.
Note that the sine basis $v_m$ is also the basis used in the Navier solution procedure
w.r.t. the bending of simply supported rectangular thin plates (cf. \cite{reddy2006theory}).

\textbf{Bi-Laplacian with all enforced boundary conditions.}
The case of all enforced boundary conditions $u=0$ and $\partial_{\nu} u  = 0$
plays a prominent role in mechanics for the problem of a clamped thin plate.
Here, the eigenvalue problem is 
\begin{align}\label{eqn:BiLapEigenAllEnf}
\laplace^2 u &= \mu u \qquad \text{in} \; \Omega, \qquad & 
u &= 0 \quad \text{and} &
\partial_{\nu} u &= 0 \quad \text{on} \; \partial \Omega.
\end{align}
We discuss this problem briefly for the square domain $\Omega$
to see that there are no separable solutions.
\begin{theorem}\label{theo:BiLapEigenAllEnf}
The eigenvalue problem stated in Eq.~\eqref{eqn:BiLapEigenAllEnf}
does not have any separable solutions.
\end{theorem}
\begin{proof}
Searching for separable solutions $u(x,y)=\alpha(x) \beta(y)$
we end up once more with Eq.~\eqref{eqn:Alpha}.
But now, the boundary conditions on $u$ at $x = -1$ tell us that
\begin{align*}
u(-1, y) &= \alpha( -1 ) \; \beta(y) = 0 \;, &
\partial_{\nu} u(-1, y) &= \alpha'( -1 ) \; \beta(y) = 0 \;,
\qquad \text{for all} \; y .
\end{align*}
Since $\beta$ must not be the zero-function for $u$ to be non-trivial, we get again to the initial value problem 
stated in Eq.~\eqref{eqn:IVP} which in turn implies $\alpha \equiv 0$ making $u$ trivial. 
Hence, this eigenproblem does not have any separable eigenfunctions at all.
\end{proof}

\textbf{Bi-Laplacian with a different set of boundary conditions,
and its relation to Hessian regularization.}
Because the regularizer $R_3$ is as used as $R_2$ in image and geometry processing,
we will also take a look at the eigenvalue problem resulting from the calculus
of variations applied to $R_3$. 
But first we compare both regularizers
\begin{align}
R_2[u] &= \frac{1}{2|\Omega|} \int\limits_\Omega \left|\laplace u \right|^2 \; dx \; , &
R_3[u] &= \frac{1}{2|\Omega|} \int\limits_\Omega \left|H u \right|_F^2 \; dx,
\end{align}
for scalar-valued functions $u$ from the space
\begin{align}\label{eqn:SpaceX2Scalar}
X_2 &:= \left\{ u:\Omega \to \R \mid \Omega = [-1,1]^2, \; u \in H^2(\Omega), \; \partial_\nu u |_{\partial \Omega} = 0\right\},
\end{align}
for which the family of tensor-product cosine functions $(u_m)_{m \in \N_0^2}$ is a basis.

\begin{theorem}\label{theo:LapReguEquivHessRegu}
On the space $X_2$ the regularizers $R_2$ and $R_3$ are the same, i.e., $R_2[u] = R_3[u]$ for any $u \in X_2$.
\end{theorem}
\begin{proof}
Given $u \in X_2$ we expand it as in Eq.~\eqref{eqn:CosExpansion} w.r.t. the cosine functions
and obtain 
\begin{align}
| H u(x) |_F^2 &= \sum\limits_{m \in \N_0^2} \sum\limits_{n \in \N_0^2} \hat{u}_m \hat{u}_n \; 
\trace \left( H u_m(x)  \; H u_n(x)  \right)\;
\end{align}
where the individual Hessian matrices also involve the sine functions of Eq.~\eqref{eqn:SinBasis} and are given by 
\begin{align}
H u_m(x)  &= -c_m \frac{\pi^2}{4}
\begin{pmatrix}
m_1^2 u_m(x) & -m_1 m_2 v_m(x) \\
-m_1 m_2 v_m(x) & m_2^2 u_m(x)
\end{pmatrix}.
\end{align}
By orthonormality of the cosine as well as the sine functions, we have
\begin{align}
\int\limits_\Omega | H u(x) |_F^2  \; dx &=
\sum\limits_{m \in \N_0^2} \left(\hat{u}_m\right)^2 \;
\int\limits_\Omega \trace \left( H u_m(x)^2 \right) \; dx.
\end{align}
Moreover, the trace expression is
\begin{align}
\trace &\left( H u_m(x)^2 \right) = \left( \partial_x^2 u_m \right)^2 + 2 \left( \partial_x \partial_y u_m \right)^2
+ \left(\partial_y^2 u_m \right)^2 \\
&= \left( \frac{\pi^2 m_1^2}{4} \right)^2 u_m^2 + 2 \frac{\pi^2 m_2^2}{4} \frac{\pi^2 m_2^2}{4} v_m^2
+ \left( \frac{\pi^2 m_2^2}{4} \right)^2 u_m^2. \nonumber
\end{align}
Since the $L^2$-norm of $u_m$ and $v_m$ is one, we get 
\begin{align}
\int\limits_\Omega \trace \left( H u_m(x)^2 \right) \; dx &= 
\left( \frac{\pi^2 m_1^2}{4} + \frac{\pi^2 m_2^2}{4} \right)^2 = \mu_m^2
\end{align}
which implies
\begin{align}
R_3[u] &= \frac{1}{2 |\Omega|} \sum\limits_{m \in \N_0^2} \mu_m^2 \;\left(\hat{u}_m\right)^2.
\end{align}
Hence, by comparison with the representation of $R_2$ in Eq.~\eqref{eqn:Reg2CosBasis},
we see that $R_2$ and $R_3$ are the same on the space $X_2$.
\end{proof}

\textbf{Note:}
this result says, that in the usual setup of image and geometry processing (i.e. $\partial_v u = 0$
on $\partial \Omega$) the regularizers $R_2$ and $R_3$ are the same when 
working on a square or rectangular domain $\Omega$.

We turn now to the eigenvalue problem based on $R_3$.
Here, the eigenvalue problem involves again the Bi-Laplacian,
but there are now different options to impose boundary conditions as in the case of $R_2$.
These can have different effects in the case of non-rectangular domains.
From the calculus of variations we obtain directly (cf. \cite{stein2018natural})
\begin{align}
\delta R_3[u] \; \varphi &= \left. \frac{d}{dh} R_3[u+h\varphi] \right|_{h=0} =
\left. \frac{1}{2|\Omega|} \frac{d}{dh} \int\limits_\Omega \left|H u + h H \varphi \right|_F^2 \; dx \right|_{h=0}
\nonumber \\
&= \frac{1}{|\Omega|} \left( \int\limits_\Omega \laplace^2 u \; \varphi \;dx 
+ \int\limits_{\partial \Omega} \nu^T \; H u \; \nabla \varphi - \partial_{\nu} \laplace u \;  \varphi \; dS(x)\right). \label{eqn:BiLapBI2}
\end{align}

Here, the first summand in the boundary integral of Eq.~\eqref{eqn:BiLapBI2} 
involves the gradient $\nabla \varphi$ and not just the directional derivative $\partial_{\nu} \varphi$.
Thus, we decompose $\nabla \varphi$ into the part $\nabla_S \varphi$ tangential
to the boundary surface $S = \partial \Omega$ and the part normal to it, i.e.,
$\nabla \varphi = \nabla_S \varphi + \partial_{\nu} \varphi \; \nu.$
With that, the boundary integral reads as
\begin{align}
\int\limits_{\partial \Omega} (\nu^T \; H u \; \nu) \partial_{\nu} \varphi + \nu^T \; H u \; \nabla_S \varphi - \partial_{\nu} \laplace u \;  \varphi \; dS(x).
\end{align}

In the next step, we recast
the boundary integral such that it involves only $\varphi$ and $\partial_{\nu} \varphi$.
Using the equivalent expression based on the surface-intrinsic divergence operator
\begin{align}\label{eqn:Surf}
\nu^T \; H u \; \nabla_S \varphi &= \diver_S \left(\varphi \; H u \; \nu \right) - \diver_S \left( \; H u \; \nu \right) \; \varphi
\end{align} 
we employ the surface-intrinsic divergence theorem to the first term on the right-hand side of Eq.~\eqref{eqn:Surf}.
But, because of $\partial S = \partial \partial \Omega = \emptyset $, the
corresponding integral term vanishes (cf. \cite{pigola2014global}).
After that, the final form of the boundary term is
\begin{align}
\int\limits_{\partial \Omega} (\nu^T \; H u \; \nu) \partial_{\nu} \varphi - \left( \diver_S \left( \; H u \; \nu \right) + \partial_{\nu} \laplace u \right) \;  \varphi \; dS(x).
\end{align}

We are again particularly interested in the case of the enforced boundary condition $\partial_{\nu} u=0$
which implies also $\partial_{\nu} \varphi=0$.
Therefore we obtain the remaining natural boundary condition given by
\begin{align}
\diver_S \left( \; H u \; \nu \right) + \partial_{\nu} \laplace u  &= 0 
\quad \text{on} \quad \partial \Omega.
\end{align}
We note that the surface-intrinsic divergence operator can be expressed in terms of an orthonormal basis 
of the tangent space of $\partial \Omega$. 
In $\R^2$ the boundary is one-dimensional. If $\tau$ is the normalized tangent vector s.t. the system
$ \{\tau, \nu\} $ is positively oriented, then we have (see also \cite{stein2018natural})
$\diver_S \left( \; H u \; \nu \right) = \tau^T \partial_{\tau} \left( \; H u \; \nu \right).$
While in $\R^3$ the boundary is two-dimensional.
If $ \{\tau_1, \tau_2 \} $ is then an orthonormal basis 
of the tangent space of $\partial \Omega$
s.t. the system $ \{\tau_1, \tau_2, \nu\} $ is positively oriented, then we have
$\diver_S \left( \; H u \; \nu \right) = \tau_1^T \partial_{\tau_1} \left( \; H u \; \nu \right)
+ \tau_2^T \partial_{\tau_1} \left( \; H u \; \nu \right).$
We consider again the square domain $\Omega = [-1,1]^2$,
then the eigenvalue problem looks different from Eq.~\eqref{eqn:EigBiLapR2} (based on $R_2$)
and reads as
\begin{align}\label{eqn:EigBiLapR3}
\laplace^2 u &= \mu u \qquad \text{in} \; \Omega, \qquad \\
\partial_{\nu} u &= 0 \qquad \text{and} \quad
\tau^T \partial_{\tau} \left( \; H u \; \nu \right) + \partial_{\nu} \laplace u = 0 \quad \text{on} \; \partial \Omega.
\nonumber
\end{align}

\begin{theorem}\label{theo:R2R3sameEig}
The eigenvalue problem given in Eq.~\eqref{eqn:EigBiLapR3} based on $R_3$
has the same solutions as the eigenvalue problem given in Eq.~\eqref{eqn:EigBiLapR2} based on $R_2$, namely the cosine basis of Eq.~\eqref{eqn:CosBasis}.
\end{theorem}
\begin{proof}
We know already that the cosine basis satisfies the PDE in Eq.~\eqref{eqn:EigBiLapR3},
the first boundary condition $\partial_{\nu} u=0$. 
Moreover, we know that the cosine basis satisfies $\partial_{\nu} \laplace u = 0$ on the boundary.
Thus, we only have to argue that the new boundary term also vanishes.
First, we consider the left boundary of the square $\Omega$, where $x = -1$.
Here, we have $\nu = -e_1$ the negative $x$-direction and $\tau = -e_2$ the positive $y$-direction. 
Using this, we obtain 
\begin{align}
\tau^T \partial_{\tau} \left( \; H u_m \; \nu \right) &= -e_2^T \partial_y (Hu_m \; e_1) = -\partial_x \partial_y^2 u_m 
\end{align}
and in terms of the basis function
\begin{align}
\tau^T \partial_{\tau} \left( \; H u_m \; \nu \right) 
& = -c_m \frac{\pi m_1}{2} \left(\frac{\pi m_2}{2}\right)^2 \sin \left( \frac{\pi m_1}{2} (x+1)\right)
 \cos \left( \frac{\pi m_2}{2} (y+1)\right).
\end{align}
Since $x$ is evaluated  at $x = -1$, we have $\tau^T \partial_{\tau} \left( \; H u_m \; \nu \right) = 0$.
An analogous argumentation gives the same result for the other three parts of the boundary.
We see that the family $(u_m)_{m \in \N_0^2}$ solves the eigenvalue problem specified in Eq.~\eqref{eqn:EigBiLapR3} 
with associated eigenvalues $\mu_m^2$.
\end{proof}

\textbf{Note:} 
for this result it is important that 
$\Omega$ is a rectangle.
In the case of a non-rectangular domain the eigenvalue 
problem~\eqref{eqn:EigBiLapR3} derived from $R_3$ 
may have different solutions than the  
eigenvalue problem~\eqref{eqn:EigBiLapR2} derived from $R_2$.

\section{Experiments}\label{sect:Experiments}

We conduct numerical experiment to support our theoretical results. 
In particular we have mentioned that the preliminary algorithm of 
Sect.~\ref{sect:PrelimAlgo} has limitations when the sampling of 
the scanning trajectory is spatially sparse, as it is for example in 
the OpenMPI Dataset \cite{knopp2020openmpidata}. 
We show in Exp.~1 that, in the case of sparse scans, 
employing the proposed second order regularization
yields consistently reconstructions of higher quality, 
when compared to the first order regularizer; 
the quality of the results is proven both qualitatively 
and quantitatively. Additionally, in Exp.~2, 
we verify that the proposed second order regularization  
is an effective quality enhancing method also with denser scans. 
Together, these two experiments demonstrate the potential 
of the proposed second order regularization.

\textbf{Parameter Selection.} 
The strength of the regularization in the core stage of the 
algorithm (Eq.~\eqref{eqn:EnergyBasic}) is controlled by 
the parameter $\lambda$. For the choice of this parameter 
we created a dataset of simulated scans and performed 
a reconstruction of each phantom for a variety of choices of $\lambda$.
Then, we computed the PSNR between the outputs $u$ 
and the ground truths $u_{\mathrm{GT}}$.
Given a ground truth phantom $\rho_{\mathrm{GT}}$
we used the ideal trace $u_{\mathrm{GT}}:=\kappa_{h}*\rho_{\mathrm{GT}}$
as ground truth for the core response.  
Finally, we selected the parameter $\lambda^*$ that maximizes the average PSNR. 
For the deconvolution stage we proceeded analogously for the regularization strength 
parameter $\mu$ in Eq.~\eqref{eqn:FunctionalJtraceGeneric}. 
More precisely, we keep $\lambda=\lambda^*$ fixed 
to obtain the input $u$ and then deconvolve 
for a range of different parameters $\mu$. Then, we compute the 
PSNR between the reconstructions and $\rho_{\mathrm{GT}}$, average 
along the dataset and display results for the parameter $\mu^*$ that 
maximizes the average PSNR computed.
The parameter $\lambda$ has been selected proceeding in two steps:  
first, we have checked the correct magnitude by performing reconstruction 
with $\lambda = j\cdot 10^{i}$ for $i\in\lbrace -3,\dots ,3\rbrace$ 
and $j\in\lbrace 1,5\rbrace$; once the correct magnitude parameter $i^*$ 
has been found, we have run another grid search around $i^*$ by 
considering $\lambda = j*10^{i}$ where $j\in\lbrace 1,\dots ,9\rbrace$ 
and $i\in\lbrace (i^*-1 )\, , i^* , (i^* +1) \rbrace$. 
For the parameter $\mu$ in the deconvolution stage, we have 
considered a grid search for the values $\mu = j\cdot 10^{i}$ 
where $j\in\lbrace 1,5\rbrace$ and $i\in\lbrace -4, 2\rbrace$.

The simulated signal data is obtained from a dataset 
consisting of phantoms that are shaped as the characters 
of an open source font package. In particular, we have 
produced 62 phantoms using all letters of the English alphabet, 
both capitalized and non capitalized, plus the digits from 0 to 9. 
The ground truths are convolved discretely on a 
$1000 \times 1000$ grid with the matrix Kernel $K_h$. 
The resolution parameter $h$ is set to $h=0.01$, which is 
within the typical regime as discussed in \cite{marz2016model}. 
The convolution result $A$ is then used according 
to Eq.~\eqref{eqn:Split} to obtain a signal $s(t)$ 
on the Lissajous curve
\begin{align}
r(t) &= \left( \sin\left(2\pi 16 t+ \frac{\pi}{2}\right)\; , \; \sin\left(2\pi 17 t+ \frac{\pi}{2}\right) \right).
\end{align}
Then, $L=1632$ equidistant time sample $t_l \in [0,1]$ 
are taken to produce the clean signal time series $s_l = s(t_l)$. 
The choice of the trajectory and its time sampling are 
inspired by the OpenMPI 2D scanning sequence, the corresponding 
curve is shown in Fig.~\ref{fig:liss:scans} (left). 
For the experiment 2, we have merged the standard scan with a 
rotated scan, to enhance to spatial sampling as described 
in \cite{gapyak2022mdpi}. The resulting sampling is 
plotted in Fig.~\ref{fig:liss:scans} (right). 
Finally, to account for noise which is always present 
in real scans, we add Gaussian noise with a noise level of 2\%. 
The signal data from which we want to reconstruct is thus
\begin{align}
\hat{s}_l &= s_l + \varepsilon \; N_l &
\varepsilon &= 0.1 \; \max \left\{ |s_l| : l=1:L\right\}.
\end{align}
The random variables $N_l$ are i.i.d normally distributed 
with zero mean and a standard deviation of one.

\subsection{Experiment 1: The Second Order Approach Improves Reconstructions on Sparse Data}\label{sect:Exp1}

\def\imratio{0.3}
\begin{figure}[t]
\centering
\captionsetup[subfigure]{aboveskip=0pt,belowskip=0pt}
\begin{subfigure}[t]{\imratio\linewidth}
	\includegraphics[width=\linewidth]{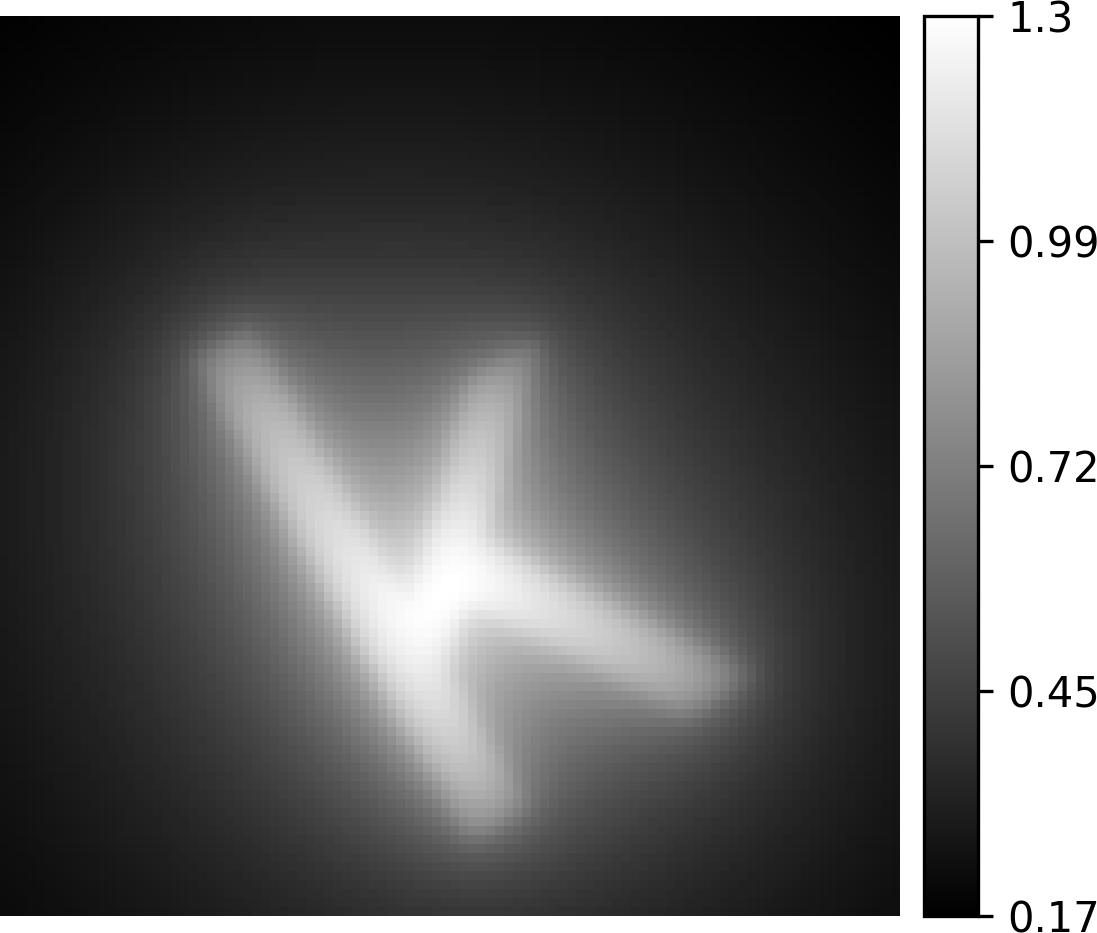}
	\caption{\centering\scriptsize Trace ground truth.}
	\label{subfig:gt:trace}
\end{subfigure}
\hfill
\begin{subfigure}[t]{\imratio\linewidth}
	\includegraphics[width=\linewidth]{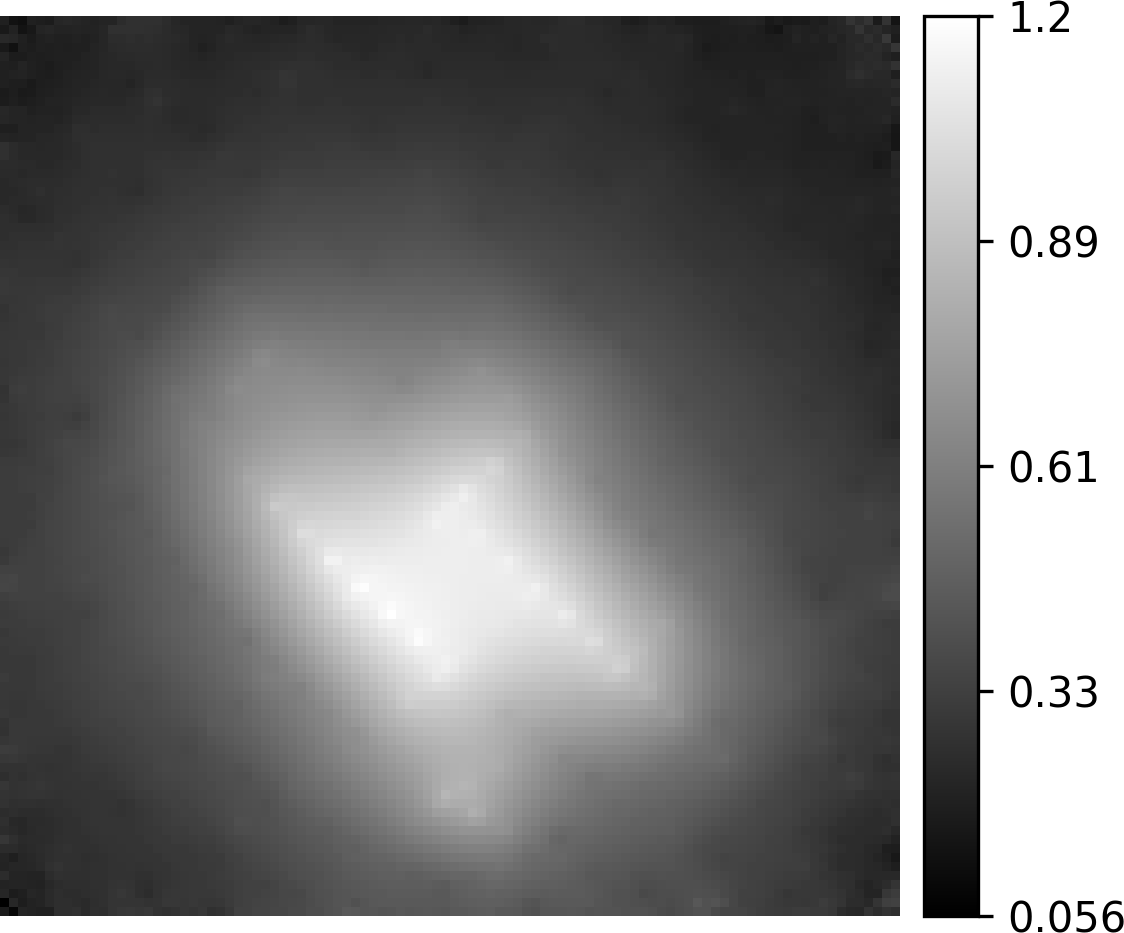}
	\caption{\centering\scriptsize Trace reco. w. 1st order. \protect\linebreak PSNR 27.5, SSIM 0.982 }
	\label{subfig:trace1sparse}
\end{subfigure}
\hfill
\begin{subfigure}[t]{\imratio\linewidth}
	\includegraphics[width=\linewidth]{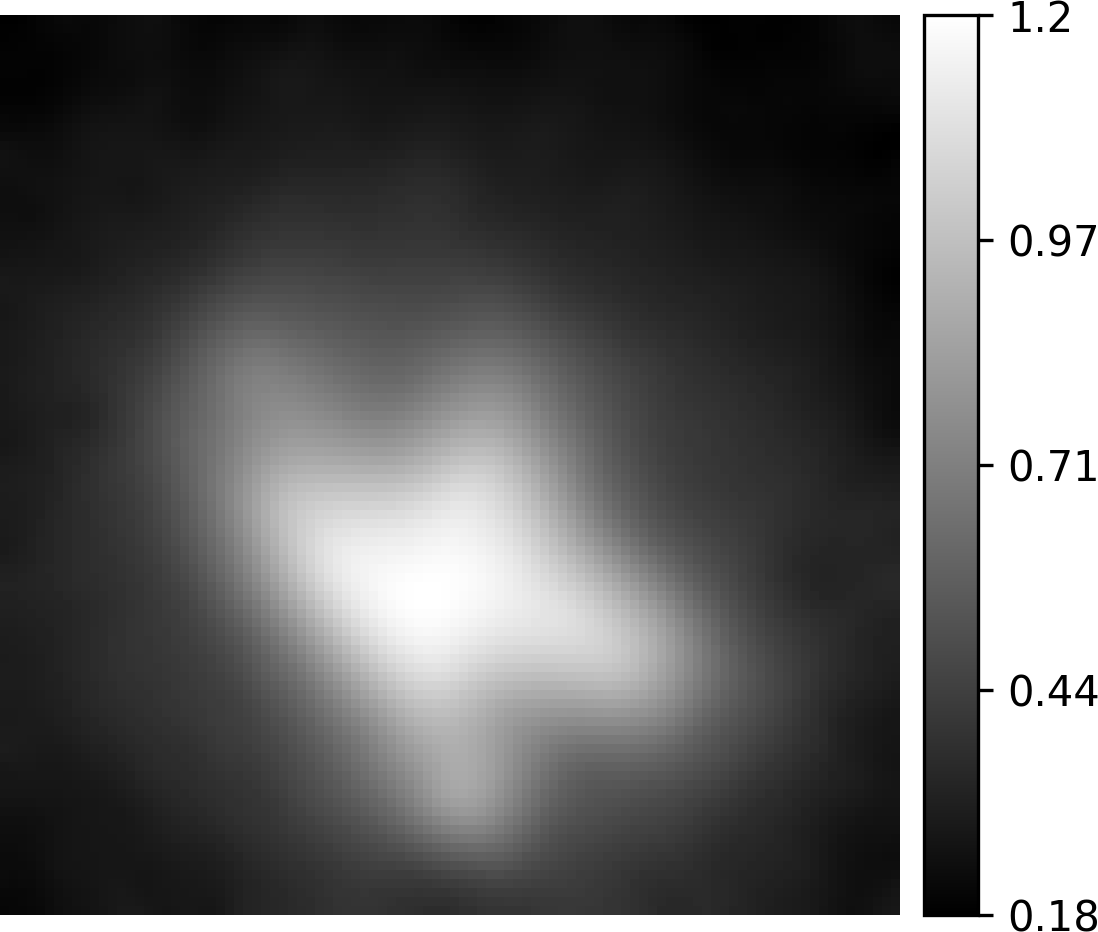}
	\caption{\centering\scriptsize Trace reco. w. 2nd order. \protect\linebreak PSNR 29.26, SSIM 0.989 }
	\label{subfig:trace2sparse}
\end{subfigure}
\begin{subfigure}[t]{\imratio\linewidth}
	\includegraphics[width=\linewidth]{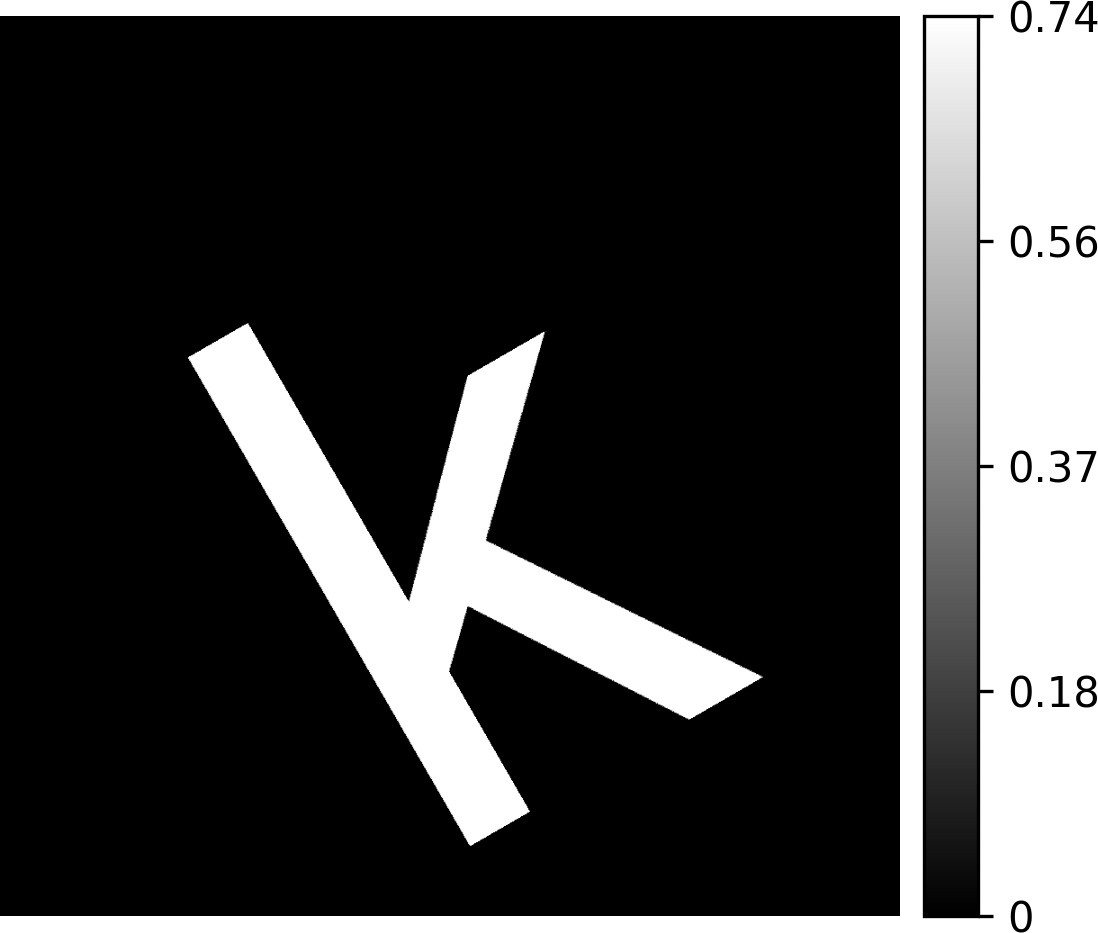}
	\caption{\centering\scriptsize Ground truth.}
	\label{subfig:gt}
\end{subfigure}
\hfill
\begin{subfigure}[t]{\imratio\linewidth}
	\includegraphics[width=\linewidth]{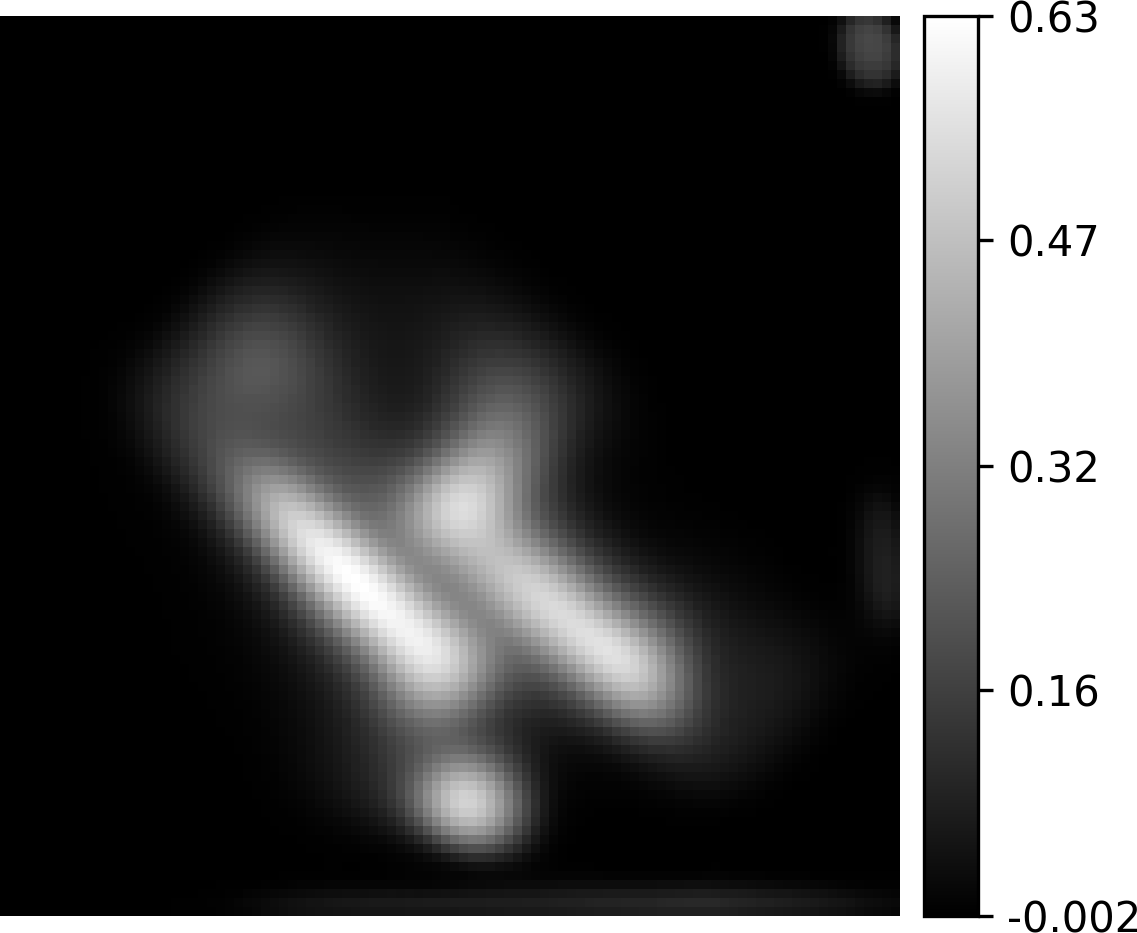}
	\caption{\centering\scriptsize Reco. from (\ref{subfig:trace1sparse}). \protect\linebreak PSNR 14.1, SSIM 0.661}
	\label{subfig:rec1sparse}
\end{subfigure}
\hfill
\begin{subfigure}[t]{\imratio\linewidth}
	\includegraphics[width=\linewidth]{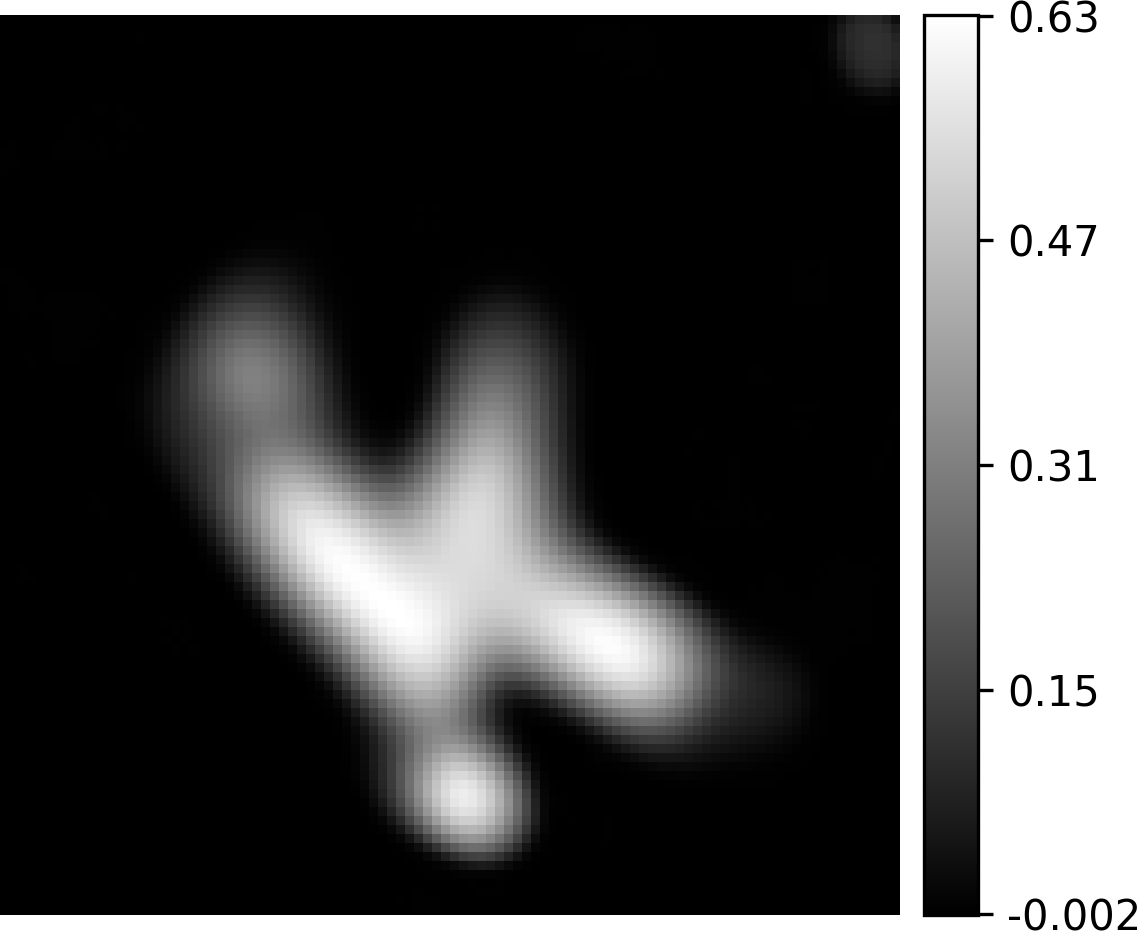}
	\caption{\centering\scriptsize Reco. from (\ref{subfig:trace2sparse}). \protect\linebreak PSNR 14.63, SSIM 0.727}
	\label{subfig:rec2sparse}
\end{subfigure}
\caption{\footnotesize 
Reconstruction of a k-shaped phantom from sparse data (Exp.~1). 
When using first order regularization, the trace is "spiky" (\ref{subfig:trace1sparse}) 
and the final reconstruction (\ref{subfig:rec1sparse}) has significant artifacts and a central disconnection. 
In contrast, second order regularization produces a smooth trace (\ref{subfig:trace2sparse}), 
resulting in a successful reconstruction (\ref{subfig:rec2sparse}) of the phantom's connected shape.
This illustrates the benefit of the proposed second order regularization.
}
\label{fig:exp1}
\end{figure}

We consider 
the case of one single Lissajous scan as shown 
in Fig.~\ref{fig:liss:scans} (left), with the aim to show the effectiveness 
of employing the second order regularizer in the case of sparse samplings,
which is a natural scenario in MPI. 
The selected parameters for the reconstruction using the first order regularizer 
in the core stage are $\lambda^* = 0.08$, and $\mu^* = 0.05$ for the deconvolution. 
For the second order regularizer in the core stage we use $\lambda^* = 0.01$, and $\mu^* = 0.01$. 
The average PSNR and SSIM scores for these parameters on the dataset are 
displayed in Tab.~\ref{tab:validation:exp1}.
We observe the relatively low absolute PSNR values of the deconvolution stage 
which is due to the specific type of images.
In fact, the ground truth data are binary images of various shapes of letters,
where about $1-2\%$ (roughly estimated) of the pixels belong to the boundary of 
the shape with jumps of jump height equal to the grayscale range. 
Compared with that, jumps in natural images have typically a fraction of the grayscale 
range as jump heights. Hence, deviations of a reconstruction near the boundary 
have a greater impact on the PSNR value (which is a kind of relative error quantity) 
resulting in lower values as usually observed for natural images.

\begin{table}[t]
	\scriptsize {
		\begin{center}
			\begin{tabular}{|c|c|c| }
				\hline
				Stage	& 	$1^{st}$ order 	& 	$2^{nd}$ order \\
				\hline
				Core Stage &
				\begingroup\renewcommand{\arraystretch}{1.3}%
				\begin{tabular}{@{}c@{}}PSNR $25.48 \pm 1.12$ \\ SSIM $0.973 \pm 0.006$ \end{tabular} 
				\endgroup	& 
				\begingroup\renewcommand{\arraystretch}{1.3}%
				\begin{tabular}{@{}c@{}}PSNR $27.26 \pm 1.07$ \\ SSIM $0.984 \pm 0.003$ \end{tabular}
				\endgroup	\\
				\hline
				Deconvolution &
				\begingroup\renewcommand{\arraystretch}{1.3}%
				\begin{tabular}{@{}c@{}}PSNR $13.71 \pm 1.380$ \\ SSIM $0.631 \pm 0.035$ \end{tabular} 
				\endgroup &
				\begingroup\renewcommand{\arraystretch}{1.3}%
				\begin{tabular}{@{}c@{}}PSNR $14.27 \pm 1.36$ \\ SSIM $0.707 \pm 0.027$ \end{tabular} 
				\endgroup \\		
				\hline
			\end{tabular}
		\end{center}
	}
	\caption{\footnotesize 
	Average reconstruction scores after the core stage (first row) and after deconvolution (second row) 
	on the dataset with the sparse sampling in Fig.~\ref{fig:liss:scans} (left). 
	The higher scores obtained with the proposed second order regularization demonstrate 
	its benefit for sparse samplings.}
	\label{tab:validation:exp1}
\end{table}

For a qualitative analysis, we display the reconstructions 
for a k-shaped phantom from the used dataset in Fig.~\ref{fig:exp1}. 
First, we look at the reconstructed traces (upper row of Fig.~\ref{fig:exp1}).
In the trace reconstructed with the first order regularizer (Fig.~\ref{subfig:trace1sparse}), 
we can see the spikes corresponding to the locations of the data points 
along the scanning trajectory. 
Moreover, we observe the empty area near the central junction point 
of the letter k; this is a consequence of the sparse sampling. 
In fact, from the spikes showing the Lissajous trajectory 
in Fig.~\ref{subfig:trace1sparse}, we can see that the junction 
point is undersampled. Comparing Fig.~\ref{subfig:trace1sparse}
with Fig.~\ref{subfig:trace2sparse},
we can see that the result of the second order approach (Fig.~\ref{subfig:trace2sparse})
is much smoother and the central junction 
area of the letter k is preserved. Moreover, when comparing 
with the ground truth trace $u_{\mathrm{GT}}$, the trace obtained with 
the second order approach in Fig.~\ref{subfig:trace2sparse} is closer 
to the ground truth in Fig.~\ref{subfig:gt:trace}; this observation is also 
supported quantitatively by the PSNR and SSIM scores of the specific example.
From the final reconstructions (lower row of Fig.~\ref{fig:exp1}), we can see 
that the reconstruction in Fig.~\ref{subfig:rec2sparse} obtained by deconvolving 
the trace in Fig.~\ref{subfig:trace2sparse} yields visually a better result 
as well as higher PSNR and SSIM 
scores, compared with the reconstruction in Fig.~\ref{subfig:rec1sparse} 
obtained from the trace in Fig.~\ref{subfig:trace1sparse}. 
This shows that 
the second order regularizer in the first stage translates to an increase 
in the overall quality of the reconstruction. The systematic improvement in 
the quality of reconstruction in the case of sparse data is quantitatively 
visible from Tab.~\ref{tab:validation:exp1}: both the average PSNR 
and SSIM of the reconstructions obtained with the second order regularizer 
are higher than those obtained using the first order regularizer. 

\subsection{Experiment 2: Improvement of the Second Order Approach on Denser Data}\label{sect:Exp2}

\def\imratio{0.30}
\begin{figure}[t]
\centering
\captionsetup[subfigure]{aboveskip=0pt,belowskip=0pt}
\begin{subfigure}[t]{\imratio\linewidth}
	\includegraphics[width=\linewidth]{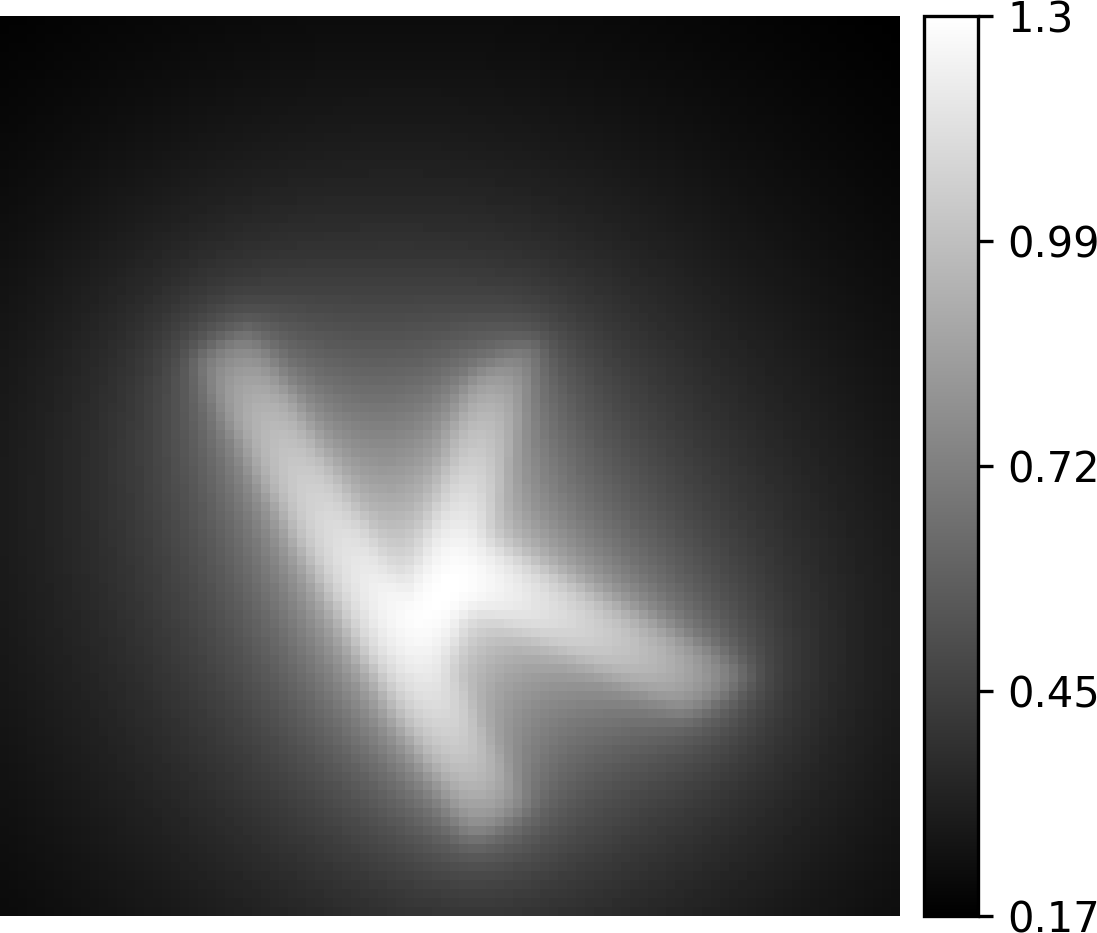}
	\caption{\centering\scriptsize Trace ground truth.}
	\label{subfig:gt:trace:dense}
\end{subfigure}
\hfill
\begin{subfigure}[t]{\imratio\linewidth}
	\includegraphics[width=\linewidth]{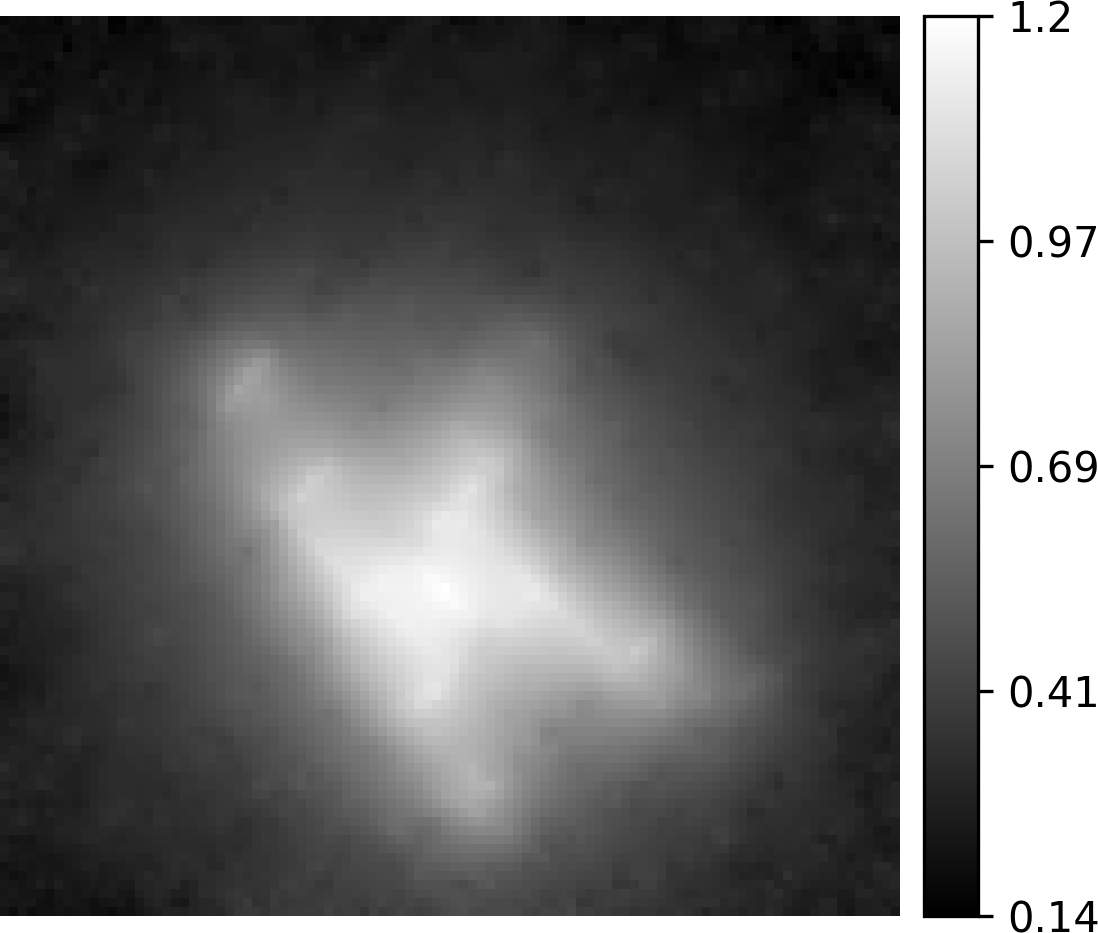}
	\caption{\centering\scriptsize Trace reco. w. 1st order. \protect\linebreak PSNR 31.46, SSIM 0.993 }
	\label{subfig:trace1dense}
\end{subfigure}
\hfill
\begin{subfigure}[t]{\imratio\linewidth}
	\includegraphics[width=\linewidth]{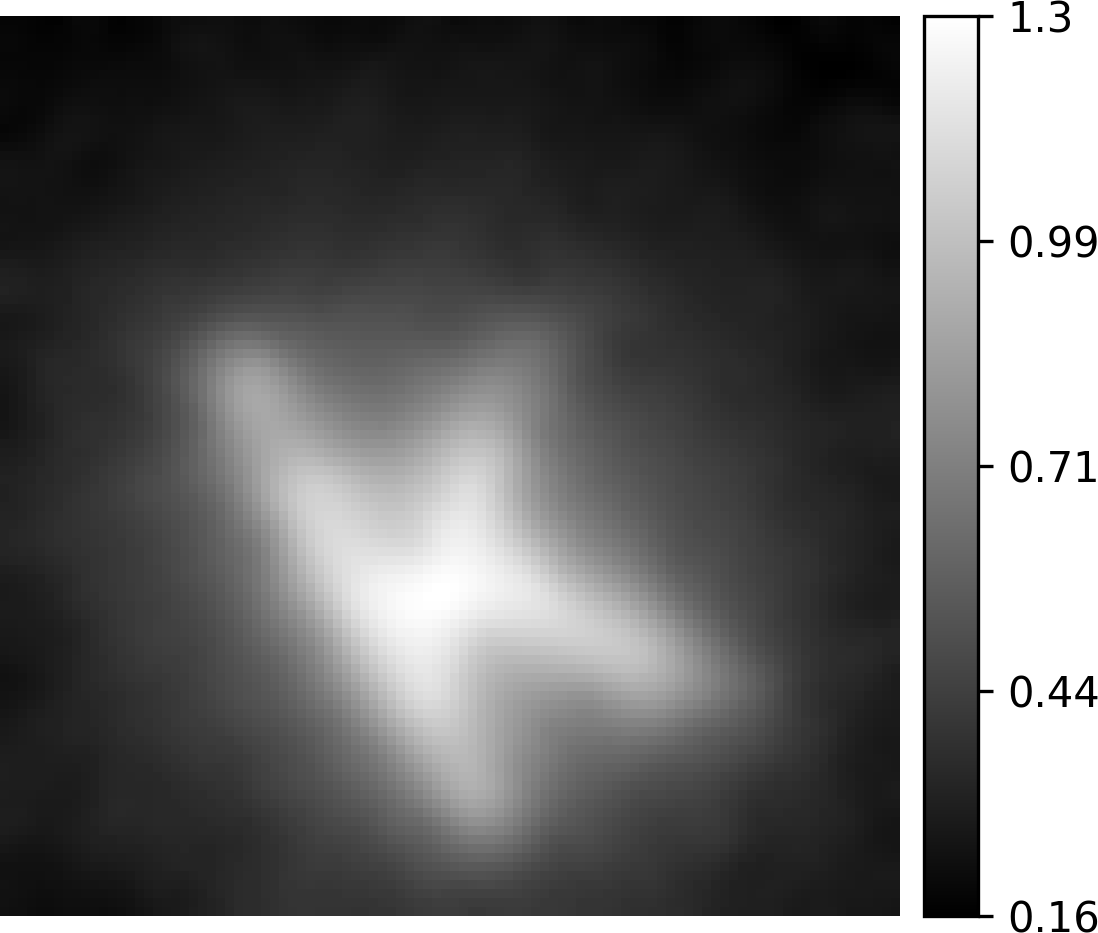}
	\caption{\centering\scriptsize Trace rec. w. 2nd order. \protect\linebreak PSNR 34.78, SSIM 0.997 }
	\label{subfig:trace2dense}
\end{subfigure}
\begin{subfigure}[t]{\imratio\linewidth}
	\includegraphics[width=\linewidth]{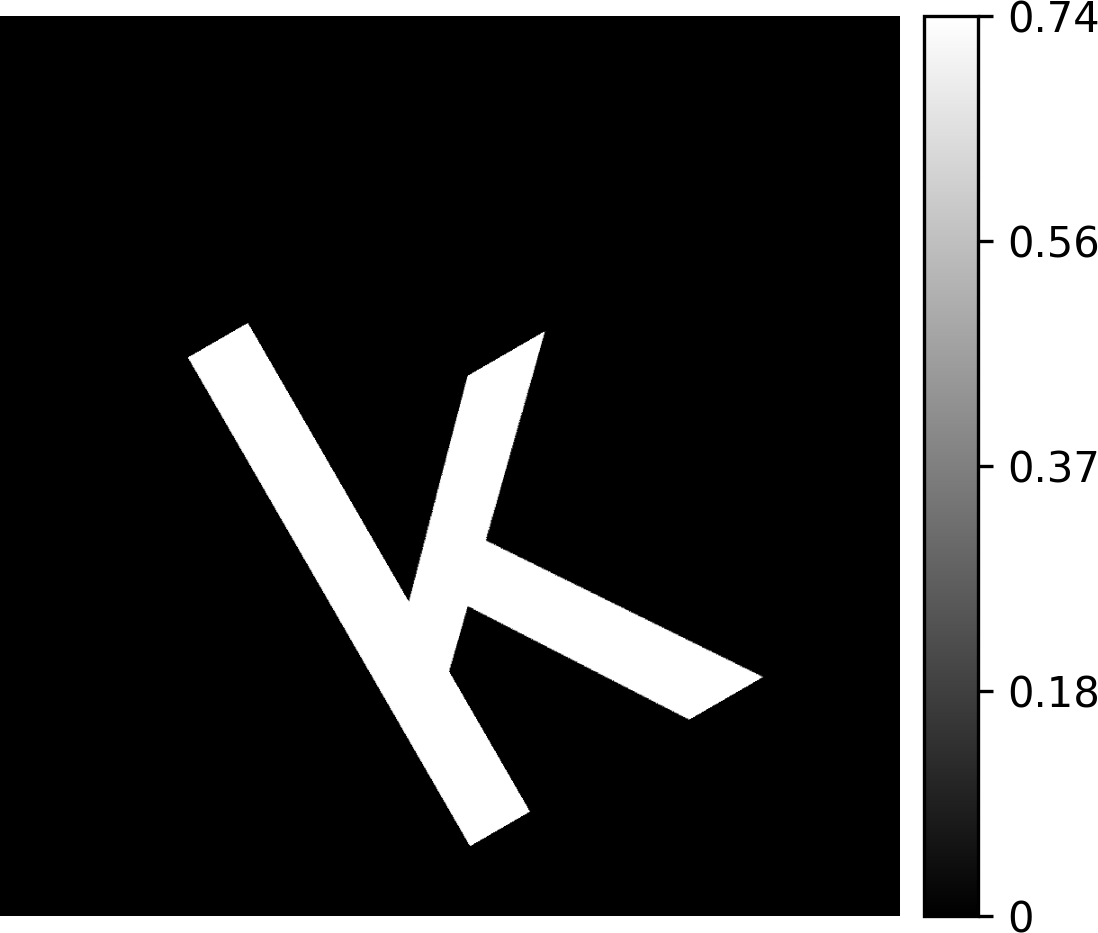}
	\caption{\centering\scriptsize Ground truth.}
	\label{subfig:gt:dense}
\end{subfigure}
\hfill
\begin{subfigure}[t]{\imratio\linewidth}
	\includegraphics[width=\linewidth]{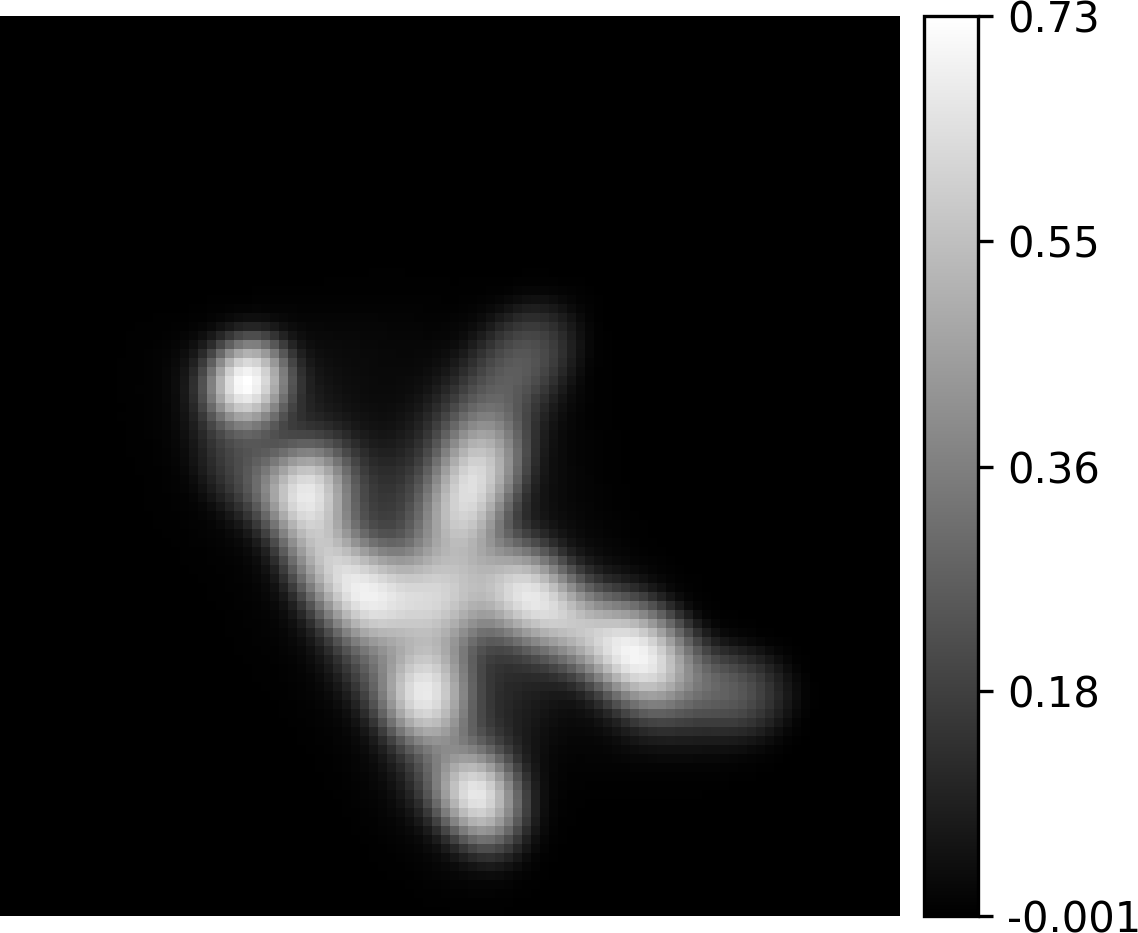}
	\caption{\centering\scriptsize Reco. from (\ref{subfig:trace1dense}). \protect\linebreak PSNR 16.31, SSIM 0.816 }
	\label{subfig:rec1dense}
\end{subfigure}
\hfill
\begin{subfigure}[t]{\imratio\linewidth}
	\includegraphics[width=\linewidth]{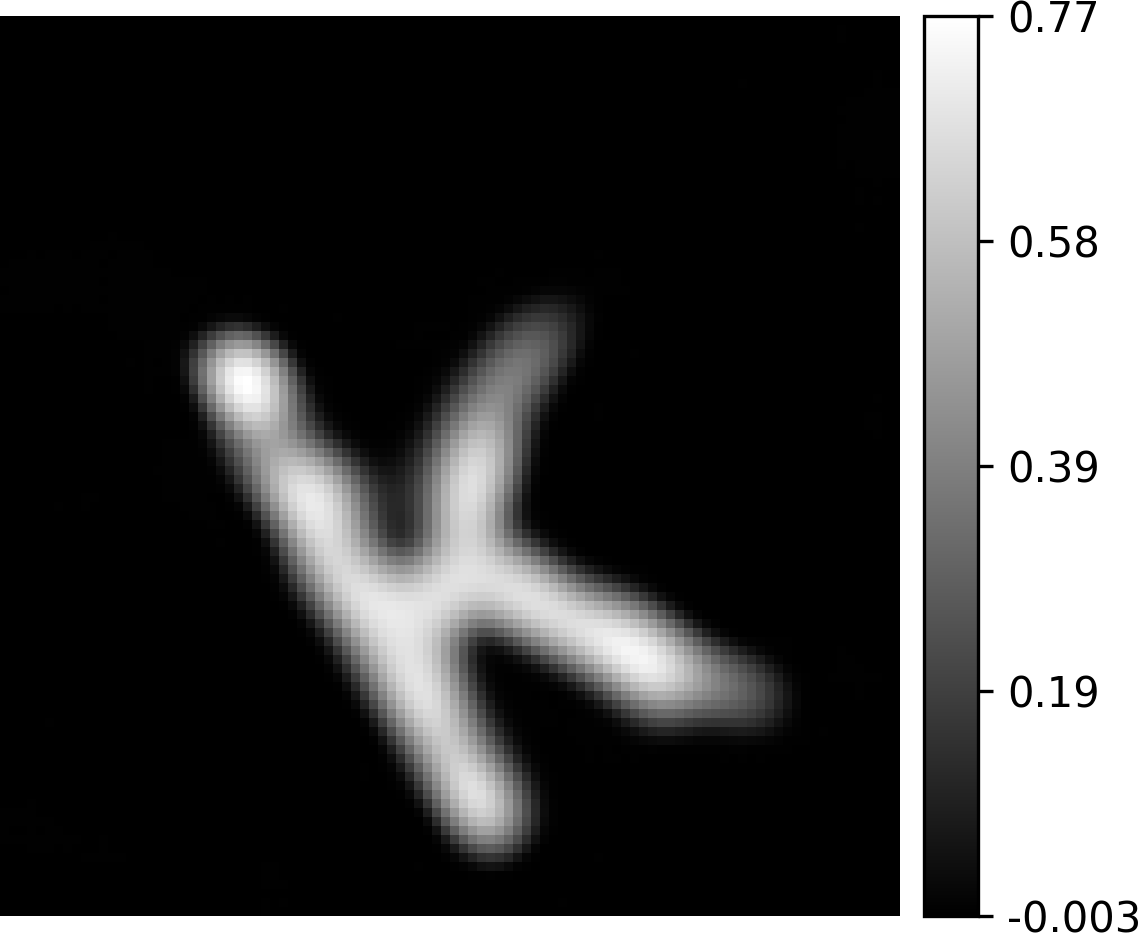}
	\caption{\centering\scriptsize Reco. from (\ref{subfig:trace2dense}). \protect\linebreak PSNR 17.79, SSIM 0.88 }
	\label{subfig:rec2dense}
\end{subfigure}
\caption{\footnotesize  
Reconstruction of a $k$-shaped phantom from denser data (Exp.~2).
The second order regularization in the core stage
yields a reconstruction (\ref{subfig:rec2dense}) with fewer 
artifacts also in the case of denser samplings, when compared with the results 
obtained using the first order regularization in the core stage (\ref{subfig:rec1dense}).}
\label{fig:exp2}
\end{figure}

\begin{table}[t]
	\scriptsize {
		\begin{center}
			\begin{tabular}{|c|c|c| }
				\hline
				Stage	& 	$1^{st}$ order 	& 	$2^{nd}$ order \\
				\hline
				Core Stage &
				\begingroup\renewcommand{\arraystretch}{1.3}%
				\begin{tabular}{@{}c@{}}PSNR $30.07 \pm 0.91$ \\ SSIM $0.992 \pm 0.001$ \end{tabular} 
				\endgroup	& 
				\begingroup\renewcommand{\arraystretch}{1.3}%
				\begin{tabular}{@{}c@{}}PSNR $33.44 \pm 0.88$ \\ SSIM $0.996 \pm 0.001$ \end{tabular}
				\endgroup	\\
				\hline
				Deconvolution &
				\begingroup\renewcommand{\arraystretch}{1.3}%
				\begin{tabular}{@{}c@{}}PSNR $16.18 \pm 1.350$ \\ SSIM $0.815 \pm 0.013$ \end{tabular} 
				\endgroup &
				\begingroup\renewcommand{\arraystretch}{1.3}%
				\begin{tabular}{@{}c@{}}PSNR $17.66 \pm 1.36$ \\ SSIM $0.880 \pm 0.014$ \end{tabular} 
				\endgroup \\		
				\hline
			\end{tabular}
		\end{center}
	}
	\caption{\footnotesize Average reconstruction score after 
	the core stage (first row) and after deconvolution (second row) 
	on the dataset with the denser sampling in Fig.~\ref{fig:liss:scans} (right). 
	The higher scores obtained with the proposed second order regularization demonstrate 
	its benefit also for more general samplings, not only for sparse ones.}
	\label{tab:validation:exp2}
\end{table}

In Exp.~1, we have shown the benefit of the proposed second order regularization
in case of sparse data. 
The proposed algorithm, however, considers the input data as a set of data points 
and can merge multiple scans to have effectively a denser sampling \cite{gapyak2022mdpi,gapyak2023multipatch}: 
in this experiment we merge the standard scan with 
the scan of the phantom rotated by $90^\circ$. 
The considered sampling in this experiment is shown in 
Fig.~\ref{fig:liss:scans} (right). In this way, we double the 
amount of points and increase the spatial density of the 
samples as suggested in \cite{gapyak2022mdpi}. 
The aim of this 
experiment is to show that second 
order regularization is beneficial also if the scans are denser,
and thus, constitutes a more general reconstruction-enhancing 
technique. To this aim, we first select the regularization 
parameter on our custom dataset: for the first order 
regularizer in the core stage we have 
selected $\lambda^* = 1$ and $\mu^* = 0.05$ for the 
deconvolution. For the second order regularizer in the 
core stage we have selected $\lambda^* = 0.004$ 
and $\mu^* = 0.01$. The average PSNR and SSIM scores are 
displayed in Tab.~\ref{tab:validation:exp2}. For visual 
assessment of the reconstruction, we have displayed the 
reconstructions obtained on the same k-shaped phantom of 
Exp.~2 in Fig.~\ref{fig:exp2}. From a qualitative 
viewpoint, the employment of the second order regularizer 
yields a trace (Fig.~\ref{subfig:trace2dense}) which 
strongly resembles the ground truth $u_{\mathrm{GT}}$ 
in Fig.~\ref{subfig:gt:trace:dense}. Moreover, the second 
order regularization yields a smoother trace when 
compared with the result obtained with the first order 
regularizer in Fig.~\ref{subfig:trace1dense}, where we can 
still see the spikes in the sampling locations and the 
undersampled areas towards the center of the FoV. 
The undersampled region and the first order regularization 
translate into a reconstruction (Fig.~\ref{subfig:rec1dense}) 
which looks disconnected in the top part of the $k$ phantom; 
the same area is successfully recovered both in the 
trace and the final reconstruction by employing the second 
order regularization (cf. Fig.~\ref{subfig:rec2dense}). 
Quantitatively, the improvement in the reconstruction 
quality by employing the second order regularization is 
also visible from the higher average PSNR and SSIM scores
computed both after the core stage and the deconvolution 
in Tab.~\ref{tab:validation:exp2}.
(Regarding the relatively small absolute PSNR values after the deconvolution stage
we refer to the discussion at the end of the first paragraph of Sect.~\ref{sect:Exp1}.)

\section{Conclusions}\label{sect:Conclusion}

In this paper we have proposed a second order Laplacian-regularization 
technique to enhance the overall reconstruction quality of 
model-based MPI-reconstruction. 
For the reconstruction of the MPI-core response
(in the core stage of our two-stage algorithm)
we have incorporated the $L^2$-norm of the Laplacian as regularizer into the energy functional
which features a model-based data fidelity term. 
We have argued that the second order regularization serves two porposes:
(i) better gap filling in the case of sparse data 
and (ii) higher smoothness of the resonconstructed MPI core response, which is an 
analytic function in theory \cite{marz2016model}.
The design of the core stage in the novel two-stage algorithm
uses tensor-product cosine functions opposed to finite difference schemes.
For a square domain the cosine functions are both the eigenfunctions 
of the neagtive Laplacian with natural boundary conditions
as well as the eigenfunctions of the Bi-Laplacian with enforced zero 
Neumann boundary conditions.
The eigenfunctions gave a simpler description for both regularizers
the $L^2$-norm of the gradient (Eq.~\eqref{eqn:RegularizerCont1}) and the
$L^2$-norm of the Laplacian (Eq.~\eqref{eqn:RegularizerR2}).
We have provided a theoretical foundation of our approach. 
In particular, we have studied in depth the eigenvalue problem for the Bi-Laplacian
on a square w.r.t. the existence of tensor product solutions and have revealed 
the importance of the imposed boundary conditions.
The choice of boundary conditions is more intricate (compared with the Laplacian case),
and it affects the eigenfunctions.  
We have shown that in the case of all natural boundary conditions 
only for the kernel of the Bi-Laplacian tensor-product eigenfunctions are available (Thm.~\ref{theo:EigBiLapAllNatSep}), 
and that the kernel contains all harmonic functions (Thm.~\ref{theo:BiLapAllNatKernel}).
Thus, enforcing at least partially boundary conditions is important.
By partially enforcing Neumann-zero or Dirchlet-zero boundary conditions,
we have obtained either cosine or sine tensor-product eigenfunctions.
In addition, we have shown that with all enforced boundary conditions, there a no tensor-product eigenfunctions 
(Thm.~\ref{theo:BiLapEigenAllEnf}) at all.
Because Hessian regularization (Eq.~\eqref{eqn:Regularizer3}) 
is a common alternative to Laplacian regularization, we studied the corresponding
eigenvalue problem of the underlying linear operator
which is again the Bi-Laplacian but with different boundary conditions.
We have shown that in our setup Laplacian and Hessian regularization 
are actually identical 
(Thm.~\ref{theo:LapReguEquivHessRegu})
and come with the same set of eingenfunctions (Thm.~\ref{theo:R2R3sameEig}).
With experiments on simulated data we have demonstrated 
that the overall reconstruction quality is in fact enhanced both
visually and in terms of metrics such as PSNR and SSIM.
Finally, we note that the proposed method was successfully applied 
to real data~\cite{gapyak2025fast}.

\section*{Acknowledgments}
We acknowledge the support of our research by the Hessian Ministry of Higher Education, Research, Science and the Arts 
within the ``Programm zum Aufbau eines akademischen Mittelbaus an hessischen Hochschulen".
A.W. acknowledges support from the German Science Foundation (DFG) within the projects INST 168/3-1, INST 168/4-1.

\bibliographystyle{siamplain}
\bibliography{references}
\end{document}